\documentstyle[11pt,leqno,amscd,amssymb,pstricks,verbatim]{amsart}
\newtheorem{thm}{Theorem}[section]
\newtheorem{lem}[thm]{Lemma}
\newtheorem{cor}[thm]{Corollary}
\newtheorem{prop}[thm]{Proposition}
\newtheorem{con}[thm]{Conjecture}

\theoremstyle{definition}
\newtheorem{example}[thm]{Example}

\newtheorem{defn}[thm]{Definition}

\newtheorem{rems}[thm]{Remarks}

\numberwithin{equation}{thm}

\def\sH{{\mathcal H}}

\def\sL{{\mathcal L}}

\def\sV{{\mathcal V}}

\def\bX{{\bold X}}

\def\BA{{\Bbb A}}
\def\BB{{\Bbb B}}
\def\BC{{\Bbb C}}

\def\BE{{\Bbb E}}
\def\BF{{\Bbb F}}
\def\BG{{\Bbb G}}
\def\BH{{\Bbb H}}
\def\BI{{\Bbb I}}

\def\BL{{\Bbb L}}
\def\BM{{\Bbb M}}
\def\BN{{\Bbb N}}

\def\BP{{\Bbb P}}
\def\BQ{{\Bbb Q}}
\def\BR{{\Bbb R}}
\def\BS{{\Bbb S}}
\def\BT{{\Bbb T}}
\def\BU{{\Bbb U}}
\def\BV{{\Bbb V}}
\def\BW{{\Bbb W}}

\def\BZ{{\Bbb Z}}
\def\Ga{\Gamma}
\def\al{\alpha}

\def\ga{\gamma}
\def\ph{\varphi}

\def\De{\Delta}
\def\la{\lambda}
\def\th{\theta}

\def\si{\sigma}
\def\le{\leqslant}
\def\ge{\geqslant}
\def\Ob{\text{\,\rm Ob}}

\def\ggp#1#2{\left[\!\!\left[\!\begin{matrix} #1\cr #2\cr\end{matrix}\!\right]\!\!\right]}

\def\BA{{\kern.5pt\Bbb A\kern.5pt}}
\def\BB{{\kern.5pt\Bbb B\kern.5pt}}
\def\BC{{\kern.5pt\Bbb C\kern.5pt}}
\def\BE{{\kern.5pt\Bbb E\kern.5pt}}
\def\BF{{\kern.5pt\Bbb F\kern.5pt}}
\def\BG{{\kern.5pt\Bbb G\kern.5pt}}
\def\BH{{\kern.5pt\Bbb H\kern.5pt}}
\def\BI{{\kern.5pt\Bbb I\kern.5pt}}
\def\BL{{\kern.5pt\Bbb L\kern.5pt}}
\def\BM{{\kern.5pt\Bbb M\kern.5pt}}
\def\BN{{\kern.5pt\Bbb N\kern.5pt}}
\def\BP{{\kern.5pt\Bbb P\kern.5pt}}
\def\BQ{{\kern.5pt\Bbb Q\kern.5pt}}
\def\BR{{\kern.5pt\Bbb R\kern.5pt}}
\def\BS{{\kern.5pt\Bbb S\kern.5pt}}
\def\BT{{\kern.5pt\Bbb T\kern.5pt}}
\def\BU{{\kern.5pt\Bbb U\kern.5pt}}
\def\BV{{\kern.5pt\Bbb V\kern.5pt}}
\def\BW{{\kern.5pt\Bbb W\kern.5pt}}
\def\BZ{{\kern.5pt\Bbb Z\kern.5pt}}
\def\GL{{\kern.5pt\Bbb{GL}\kern.5pt}}
\def\SL{{\kern.5pt\Bbb{SL}\kern.5pt}}

\def\Ga{\Gamma}
\def\al{\alpha}

\def\ga{\gamma}
\def\ph{\varphi}

\def\De{\Delta}
\def\la{\lambda}
\def\th{\theta}

\def\si{\sigma}
\def\EE#1#2#3#4#5#6{{\!\raise1.5pt\hbox{$\begin{matrix}{\scriptstyle #1#2}\hfill\cr
\noalign{\vskip-8pt}{\scriptstyle #3#4}\hfill\cr
\noalign{\vskip-8pt}{\scriptstyle #5#6}\hfill\end{matrix}$}}}

\def\End{\operatorname{End}}

\def\Hom{\operatorname{Hom}}
\def\Tr{\operatorname{Tr}}

\def\Ker{\operatorname{Ker}}

\def\modh{\text{{\bf mod}-}}

\def\modf{\text{{\bf mod$^F$}-}}

\def\op{{\operatorname{op}}}

\newfam\eusmfam
\newfam\cmmibfam
\def\Catname#1#2{{\fam\eusmfam\relax#1#2}
}

\def\RR{{\Catname K{}\kern.5pt}}
\def\kk{{\Catname k{}\kern.5pt}}

\def\Ob{\operatorname{Ob}}
\def\Hom{\operatorname{Hom}}

\def\Ker{\operatorname{Ker}}
\def\End{\operatorname{End}}

\def\op{{\operatorname{op}}}

\def\Aut{\operatorname{Aut}}

\def\Tr{\operatorname{Tr}}

\def\ggp#1#2{\left[\kern-3.2pt\left[{#1\atop #2}\right]\kern-3.2pt\right]}
\def\sggp#1#2{\left[\kern-2.5pt\left[{#1\atop #2}\right]\kern-2.5pt\right]}
\def\ggi#1{\left[\kern-1.8pt\left[#1\right]\kern-1.8pt\right]}

\def\Catname#1#2{{\fam\eusmfam\relax#1#2}
}

\def\ph{\varphi}

\def\[[{{[\![}}
\def\]]{{]\!]}}

\def\Irr{{\text{\rm Irr}}}

\def\diag{\text{\rm diag}}

\def\kk{{\Catname {}{k}}}

\def\Irr{{\text{\rm Irr}}}

\def\diag{\text{\rm diag}}

\def\kk{{\Catname {}{k}}}

\def\lr#1{\langle #1\rangle}
\def\az{\alpha}  \def\thz{\theta}
\def\bz{\beta}  \def\dt{\Delta}
\def\gz{\gamma}
\def\ggz{\Gamma}  \def\ooz{\Omega}
\def\sz{\sigma}  \def\hsz{{\hat \sz}}

\def\ez{\varepsilon}

\def\lz{\lambda}
\def\ct{{\cal T}}  \def\fs{{\frak s}}    
\def\cc{{\cal C}}
\def\cq{{\cal Q}}     
 \def\co{{\cal O}} \def\cp{{\cal P}}
  \def\cd{{\cal D}}
\def\ca{{\cal A}}  \def\ci{{\cal I}}
\def\bbn{{\mathbb N}}  \def\bbz{{\mathbb Z}}  \def\bbq{{\mathbb Q}}
    \def\bbf{{\mathbb F}}
     
 \def\vphi{\varphi}
        \def\bd{{\bf d}} 
\def\bi{{\bf i}} \def\bj{{\bf j}}
  \def\bk{{\bf k}}

\def\brho{{\boldsymbol \rho}}
\def\ofq{{\overline{\bbf}_q}}
\def\fqi{{\bbf_{q^{\ez_\bi}}}}   \def\fqj{{\bbf_{q^{\ez_\bj}}}}
\def\fqaz{{\bbf_{q^{\ez_\brho}}}}
  \def\leq{\leqslant}  \def\geq{\geqslant}
\def\ge{\geqslant}  \def\le{\leqslant}
\def\lra{\longrightarrow}   \def\skrel{\stackrel}
\def\lmto{\longmapsto}      \def\ol{\overline}
\def\ra{\rightarrow}
\def\hom{\mbox{\rm Hom}}
\def\op{\mbox{\scriptsize op}}
   \def\Ker{\mbox{\rm Ker}\,} \def\Irr{\mbox{\rm Irr}}
\def\dim{\mbox{\rm dim}\,}   \def\Coker{\mbox{\rm Coker}\,}
\def\udim{{\mathbf dim\,}}
\def\mod{\text{{\bf mod}-}}
\def\undermod{\underline {\text{\bf mod}}\text{-}}
\def\overmod{\overline {\text{\bf  mod}}\text{-}}
\def\undermodf{\underline {\text{\bf mod}}^F\text{-}}
\def\overmodf{\overline {\text{\bf  mod}}^F\text{-}}
\def\modf{\text{{\bf mod}$^F$-}}   \def\top\Tr{\mbox{\rm Tr}\,}
 \def\soc{\mbox{\rm Soc\,}}
\def\rad{\mbox{\rm Rad}\,}  \def\top{\mbox{\rm top}\,}
\def\rep{\mbox{\rm Rep}\,}  \def\End{\text{\rm End}}

  \def\tM{{\tilde M}}  \def\tN{{\tilde N}} \def\tL{{\tilde L}}
  \def\bX{{X_k}}
\def\fM{{\frak M}}

\begin{document}

\title[Frobenius Morphisms and Representations of algebras]%
{Frobenius morphisms and Representations of algebras}
\author{Bangming Deng and Jie Du}
\address{Department of Mathematics, Beijing Normal University,
Beijing 100875,  China.} \email{dengbm@@bnu.edu.cn}
\address{School of Mathematics, University of New South Wales,
Sydney 2052, Australia.}
 \email{j.du@@unsw.edu.au\quad{\it Home
Page: {\tt http://www.maths.unsw.edu.au/$\sim$jied}}}
\date{8 July, 2003}

\thanks{Supported partially by the NSF of China,  the TRAPOYP, and the
Australian Research Council.}

\begin{abstract}
By introducing Frobenius morphisms $F$ on algebras $A$ and their
modules over the algebraic closure ${{\overline \BF}}_q$ of the
finite field $\BF_q$ of $q$ elements, we establish a relation
between the representation theory of $A$ over ${{\overline
\BF}}_q$ and that of the $F$-fixed point algebra $A^F$ over
$\BF_q$. More precisely, we prove that the category $\modh A^F$ of
finite dimensional $A^F$-modules is equivalent to the subcategory
of finite dimensional $F$-stable $A$-modules, and, when $A$ is
finite dimensional, we establish a bijection between the
isoclasses of indecomposable $A^F$-modules and the $F$-orbits of
the isoclasses of indecomposable $A$-modules. Applying the theory
to representations of quivers with automorphisms, we show that
representations of a modulated quiver (or a species) over $\BF_q$
can be interpreted as $F$-stable representations of a
corresponding quiver over ${{\overline \BF}}_q$. We further prove
that every finite dimensional hereditary algebra over $\BF_q$ is
Morita equivalent to some $A^F$, where $A$ is the path algebra of
a quiver $Q$ over ${{\overline \BF}}_q$ and $F$ is induced from a
certain automorphism of $Q$. A close relation between the
Auslander-Reiten theories for $A$ and $A^F$ is established. In
particular, we prove that the Auslander-Reiten (modulated) quiver
of $A^F$ is obtained by ``folding" the Auslander-Reiten quiver of
$A$. Finally, by taking Frobenius fixed points, we are able to
count the number of indecomposable representations of a modulated
quiver with a given dimension vector and to establish part of
Kac's theorem for all finite dimensional hereditary algebras over
a finite field.
\end{abstract}

\maketitle

\begin{center}
CONTENTS
\end{center}
\begin{itemize}
\item[1.]Introduction
\item[2.]$\BF_q$-structures on vector spaces
\item[3.]Algebras with Frobenius morphisms
\item[4.]Twisting modules with Frobenius maps
\item[5.]$F$-periods and indecomposable $F$-stable modules
\item[6.]Finite dimensional hereditary algebras
\item[7.]Almost split sequences
\item[8.]The Auslander-Reiten quivers
\item[9.]Counting the number of $F$-stable representations
\item[10.]Roots and indecomposable $F$-stable representations
\end{itemize}

\section{Introduction}

In his work \cite{Ga}, Gabriel introduced the idea of quiver
representations and discovered a remarkable connections between
the indecomposable representations of (simply-laced) Dynkin
quivers and the positive roots of the corresponding finite
dimensional simple Lie algebras. The theory of quiver
representations not only may be viewed as a new language for the
whole range of problems in linear algebra, but also provided a
platform (over an algebraically closed field) for bringing new
ideas and techniques from algebraic geometry and Lie theory into
the subject. There are two major fundamental developments after
Gabriel's work. In order to complete Gabriel's classification to
include all Dynkin graphs, Dlab and Ringel \cite{DR1} studied
representations of a modulated quiver (or a species) and proved
that a modulated quiver admits only finitely many indecomposable
representations if and only if it is of Dynkin type. In the
subsequent works \cite{DF,N,DR2}, representations of both quivers
and modulated quivers of ``tame'' representation type are
classified. On the other hand, by employing the methods of
invariant theory, Kac \cite{K1,K2} was able to establish a
connection between the indecomposable representations of {\it any}
finite quiver and the positive roots of the corresponding
Kac-Moody algebra with a symmetric generalized Cartan matrix.

These fundamental works provide two major approaches in the
representation theory of algebras. The quiver approach is usually
used to study representations of algebras over an algebraically
closed field in which methods in algebraic geometry and invariant
theory can be applied. However, in the context of Lie algebras, it
only deals with the case of symmetric (generalized) Cartan
matrices. The modulated quiver approach, though a bit artificial,
is purely algebraic, and suitable for an {\it arbitrary} ground
field. If the underlying valued quiver of a modulated quiver is
not constantly valued, then the corresponding symmetrizable Cartan
matrix is not symmetric. Thus, this approach covers the case of
all symmetrizable Cartan matrices.

Another important achievement in representation theory of algebras
is the discovery of Auslander-Reiten sequences by Auslander and
Reiten in 1970's. Such sequences reflect the additional structures
imposed on the category of finite dimensional modules over an
algebra by the existence of kernels and cokernels. The
Auslande-Reiten theory soon became a fundamental tool in the
study of representations of algebras.

It is well-known in the Lie theory that a non-symmetric Cartan
matrix can be obtained by ``folding'' the graph of a symmetric
Cartan matrix via a graph automorphism. Such an idea has been used
in \cite{T,L93,L98,Re,Hub2} to study representations of quivers
with an automorphism. In this paper we shall extend this idea,
combining with the idea of Frobenius morphisms in the theory of
algebraic groups, to build a direct bridge between the quiver and
modulated quiver approaches. By introducing Frobenius morphisms of
algebras $A$ defined over $\ofq$, we shall prove that
representations of the fixed-point algebra $A^F$ are obtained by
taking fixed points of $F$-stable representations of $A$. In
particular, if $A$ is the path algebra of a quiver $Q$ which
admits an admissible automorphism $\si$ and $F$ is the Frobenius
morphism on $A$ induced from $\si$, then $A^F$ is isomorphic to
the tensor algebra of the modulated quiver obtained by folding $Q$
and $A$ through $\si$ and $F$, respectively. Thus, the
representation theory of $\BF_q$-modulated quiver can be realized
as that of an ordinary quiver $Q$ by simply taking fixed-points of
$F$-stable representations of $Q$ over $\ofq$. Further, we
establish a relation between the Auslander-Reiten theories of $A$
and its fixed point algebra $A^F$.

We organize the paper as follows. Section 2 is a brief
introduction on $\bbf_q$-structures of vector spaces. In \S3 we
consider algebras $A$ with Frobenius morphisms $F$ and define
$F$-stable $A$-modules. Then we show that the category of
$F$-stable $A$-modules is isomorphic to the category of modules
over the fixed point algebra $A^F$. In \S4 we define the
(Frobenius) twist of an $A$-module and introduce the notion of
$F$-periodic $A$-modules. As a result, we prove in \S5 that each
indecomposable $A^F$-module can be obtained by ``folding''
$F$-periodic $A$-modules. In particular, if $A$ is finite
dimensional, there is a bijection between indecomposable
$A^F$-modules and $F$-orbits of the indecomposable $A$-modules. As
a first application to representations of quivers with
automorphisms, it is shown in \S6 that the representation theory
of  a modulated quiver (or a species) over $\BF_q$ is completely
determined by the representation theory of the corresponding
quiver over ${{\overline \BF}}_q$. We further prove that every
finite dimensional hereditary algebra over $\BF_q$ is Morita
equivalent to some $A^F$, where $A$ is the path algebra of a
quiver $Q$ over ${{\overline \BF}}_q$ and $F$ is the Frobenius
morphism induced from an automorphism of $Q$. In \S7 we establish
a relation between almost split sequences of $A$ and $A^F$.
Section 8 is devoted to proving that the Auslander-Reiten
(modulated) quiver of $A^F$ is obtained by ``folding" the
Auslander-Reiten quiver of $A$. In \S9 and \S10, we present
formulae of the number of indecomposable $F$-stable
representations of an ad-quiver with a fixed dimension vector and
prove part of Kac's theorem for all finite dimensional hereditary
algebras over a finite field. Thus, we reobtain and generalize
some results in \cite{Hua1,DX,Hub1}.

In \cite{DD1,DD2}, a strong monomial basis property for quantized
enveloping algebras of simply-laced Dynkin or cyclic type was
discovered. The present work was motivated by seeking a similar
property in the non-simply-laced case. Though methods such as
representation varieties and generic extensions employed in
\cite{DD1,DD2} are no longer valid when working over a finite
field, we want a new approach that carries them over. In a
forthcoming paper, we shall apply the theory developed in the
paper to complete our tasks on the strong monomial basis property
for {\it all} quantized enveloping algebras of finite type.

{\it Throughout, let $q$ be a prime power, $\BF_q$ the finite
field of $q$ elements and $k$ the algebraic closure $\ofq$ of
$\BF_q$. For any $r\ge 1$, let $\bbf_{q^r}$ denote the unique
extension field of $\bbf_q$ of degree $r$ contained in
$k=\overline{\bbf}_q$. All modules considered are left modules of
finite dimension over the base field. If $M$ is a module, $[M]$
denotes the class of modules isomorphic to $M$, i.e., the isoclass
of $M$. For any field ${\tt k}$, the notation ${\tt k}^{m\times
n}$ denotes the set of all $m\times n$ matrices over $\tt k$.}

\noindent {\bf Acknowledgement.} The main results of the paper
were presented by the second author at the {\it Conference for
Representations of Algebraic Groups}, Aarhus, June 2--6, 2003, and
by both authors at the {\it Conference on Representation Theory},
Canberra, June 30--July 4, 2003. We would like to thank the
organizers for the opportunity of attending the conferences. The
second author also thanks Wilberd van der Kallen for discussions
on Frobenius morphisms over infinite dimensional vector spaces.

\section{$\BF_q$-structures on vector spaces}

An $\BF_q$-{\it structure} on a vector space $V$ over $k$ is an
$\BF_q$-subspace $V_0$ of $V$ viewed as a space over $\BF_q$ such
that the canonical homomorphism $V_0\otimes_{\BF_q}k\to V$ is an
isomorphism. We shall always identify $V$ with $V_0\otimes k$ in
the sequel.

\begin{lem}\label{frob} A $k$-space $V$ has an $\BF_q$-structure $V_0$ if and
only if
$$V_0=V^F:=\{v\in V\mid F(v)=v\}$$
 for some $\BF_q$-linear
isomorphism $F:V\to V$ satisfying
\begin{itemize}
\item[(a)]  $F(\la v)=\la^qF(v)$ for all $v\in V$ and $\la\in k$;
\item[(b)] for any $v\in V$, $F^n(v)=v$ for some $n>0$.
\end{itemize}
\end{lem}

\begin{pf} If $V$ has an $\BF_q$-structure $V_0$, then we have
$V=V_0\otimes_{\BF_q}k$ and define $F:V\to V$ by sending
$v\otimes a$ to $v\otimes a^q$. Clearly, $V_0=V^F$ and $F$
satisfies the conditions (a) and (b). The proof for the converse
is given in \cite[3.5]{DM}.
\end{pf}

The map $F$ is called a {\it Frobenius map}. By the lemma, we see
that an $\BF_q$-structure on $V$ is equivalent to the existence of
a Frobenius map. However, different Frobenius morphisms may define
the same $\BF_q$-structure. From the proof we see that if $F:V\to
V$ is a Frobenius map, then there is a basis $\{v_i\}$ of $V$ such
that $F(\sum_i\la_iv_i)=\sum_i\la_i^qv_i$.

\begin{cor} If $F$ and $F'$ are Frobenius maps on a {\rm finite dimensional}
space $V$, then $F'\circ F^{-1}$ is linear on $V$ and there is a
positive integer $n$ such that $F^n={F'}^n$.
\end{cor}
\begin{pf} Since $V=V^F\otimes_{\BF_q}k=V^{F'}\otimes_{\BF_q}k$
and $V$ is finite dimensional, there is a positive integer $n$
such that
$$V(\BF_{q^n}):=V^F\otimes_{\BF_q}\BF_{q^n}=V^{F'}\otimes_{\BF_q}\BF_{q^n}.$$
Choose bases $\{v_i\}$ and $\{w_i\}$ for $V^F$ and $V^{F'}$,
respectively, and write $w_j=\sum_jx_{ij}v_i$ with
$x_{ij}\in\BF_{q^n}$. Now, one checks easily that $F^n=F^{\prime
n}$ on $V(\BF_{q^n}$), and hence, on $V$.
\end{pf}

Let ${\sV}_{k,\BF_q}$ be the category whose objects are vector
spaces over $k$ with {\it fixed} $\BF_q$-structures and whose
morphisms are linear maps $f:V\to W$ defined over $\BF_q$, namely,
$f(V_0)\subseteq W_0$. Clearly, if $F_V$ and $F_W$ are the
Frobenius maps on $V$ and $W$ defining $V_0$ and $W_0$,
respectively, then $f$ is defined over $\BF_q$ if and only if
$F_W\circ f=f\circ F_V$. In particular, if $F$ and $F'$ are two
Frobenius maps on $V$, then any $\BF_q$-linear isomorphism from
$V^F$ to $V^{F'}$ induces a $k$-linear isomorphism $\th$ on $V$
such that $\th\circ F=F'\circ\th$. In other words, up to
isomorphism in $\sV_{k,\BF_q}$, the $\BF_q$-structure on $V$ is
unique. Note that $\sV_{k,\BF_q}$ is an abelian category.

We shall represent the fixed $\BF_q$-structure on a $k$-space $V$
by a Frobenius map $F_V$. Thus, the objects $V$ in
${\sV}_{k,\BF_q}$ are, sometimes, written as $(V,F_V)$, and
$$\Hom_{\sV_{k,\BF_q}}(V,W)=\{f\in\Hom_k(V,W)\mid F_W\circ f=f\circ
F_V\}.$$

For $(V,F_V), (W,F_W)\in\text{Ob}(\sV_{k,\BF_q})$, there is an
induced $\BF_q$-linear map
$$F_{(V,W)}:\Hom_k(V,W)\to\Hom_k(V,W);f\mapsto
F_{(V,W)}(f)=F_W\circ f\circ F_V^{-1}.$$
 Note that, when both $V$ and $W$ are {\it infinite dimensional},
$F_{(V,W)}$ is not necessarily a Frobenius map on $\Hom_k(V,W)$ in
the sense of Lemma \ref{frob}. However, we have the following.

\begin{lem} \label{Lang}
For $(V_1,F_1), (V_2,F_2)\in\text{\rm Ob}(\sV_{k,\BF_q})$, let
$F=F_{(V_1,V_2)}$ be defined as above.

(a) We have an $\BF_q$-space isomorphism
$$
\Hom_{\BF_q}(V_1^{F_1},V_2^{F_2})\cong
\Hom_{\sV_{k,\BF_q}}(V_1,V_2)=\Hom_k(V_1,V_2)^F.$$

(b) Let $\text{\rm hom}_k(V_1,V_2)=\Hom_k(V_1,V_2)^F\otimes k.$
Then, $F$ is a Frobenius map on $\text{\rm hom}_k(V_1,V_2)$ and,
if $V_i$ are finite dimensional, then $\text{\rm
hom}_k(V_1,V_2)=\hom_k(V_1,V_2)$.
\end{lem}
\begin{pf} The isomorphism in (a) is defined by sending $f$ to
$f\otimes 1$ with inverse defined by restriction. To see (b), it
is clear that $F$ is an $\BF_q$-linear map satisfying
\ref{frob}(a) and $F$ satisfies \ref{frob}(b) on $\text{\rm
hom}_k(V_1,V_2)$. Finally, the last assertion follows from a
dimension comparison.
\end{pf}

\begin{rems} \label{Lang1}
(1) Every Frobenius map on a finite dimensional $k$-space $V$
induces the Frobenius map $F=F_{(V,V)}$ on the algebra $\End_k(V)$
of all linear transformations on $V$. Clearly, $F$ is an
$\BF_q$-algebra automorphism and induces an $\BF_q$-algebra
isomorphism $\End_k(V)^F\cong \End_{\BF_q}(V^F)$.

(2) The restriction of $F$ on $\End_k(V)$ to the general linear
group $GL(V)$ induces a Frobenius morphism on the connected
algebraic group $GL(V)$ to which the following Lang-Steinberg
theorem apply.
\end{rems}

\begin{thm} \label{lang} ({\bf Lang-Steinberg}) Let $G$ be a connected affine
algebraic group and let $F$ be a surjective endomorphism of $G$
with a finite number of fixed points. Then the map $\sL:g\mapsto
g^{-1}F(g)$ from $G$ to itself is surjective.
\end{thm}

\section{Algebras with Frobenius morphisms}

Let $A$ be a $k$-algebra with identity 1. We do not assume
generally that $A$ is finite dimensional. A map $F_A:A\to A$ is
called a {\it Frobenius morphism} on $A$ if it is a Frobenius map
on the $k$-space $A$ and it is also an $\BF_q$-algebra isomorphism
sending 1 to 1. 
For example, for a finite dimensional $k$-space $V$,
$A=\End_k(V)$ admits a Frobenius morphism $F_{(V,V)}$ induced from
a Frobenius map $F_V$ on $V$ (see \ref{Lang1}).

Given a Frobenius morphism $F_A$ on $A$, let
$$A^F:=A^{F_A}=\{a\in A\mid F_A(a)=a\}$$
be the set of $F_A$-fixed points. Then $A^F$ is an
$\BF_q$-subalgebra of $A$, and $A=A^F\otimes k$. The Frobenius
morphism $F_A$ induces an algebra isomorphism on $A^F$ and
$F_A(a\otimes \la)=F_A(a)\otimes \la^q$ for all $a\in A^F,\la\in
k$.

Let $M$ be a finite dimensional $A$-module and let
$\pi=\pi_M:A\to\End_k(M)$ be the corresponding representation. We
say that $M$ is $F_A$-{\it stable} (or simply $F$-stable) if there
is an $\BF_q$-structure $M_0$ of $M$ such that $\pi$ induces a
representation $\pi_0:A^F\to\End_{\BF_q}(M_0)$ of $A^F$. Clearly,
by Lemma \ref{frob}, $M$ is $F$-stable if and only if there is a
Frobenius map $F_M:M\to M$ such that
\begin{equation}\label{fstab}
F_M(am)=F_A(a)F_M(m),\,\,\text{ for all }a\in A,m\in M.
\end{equation}
In terms of the corresponding representation $\pi$,
the $F$-stability of $M$ simply means that there is a Frobenius
map $F_M$ on $M$ such that $\pi\circ F_A=F_{(M,M)}\circ \pi$. {\it
In the sequel, we shall fix such an $\BF_q$-structure $M_0$ for an
$F$-stable module and represent it by a Frobenius map $F_M$.}
Thus, if $M$ is an $F$-stable $A$-module with respect to $F_M$,
then $M=M^F\otimes_{\bbf_q}k$, $F_M(m\otimes\lz)=m\otimes\lz^q$
for all $m\in M^F$, $\la\in k$, and $M^F$ is an $A^F$-module.
Here, again for notational simplicity, we write $M^F$ for
$M^{F_M}$.\footnote{It should be understood that the $F$'s in
$A^F$ and $M^F$ are not the same.} We shall also use, sometimes,
the notation $(M,F_M)$ for an $F$-stable module $M$.

\begin{lem} \label{IND} Let $(M_1,F_1)$ and $(M_2,F_2)$ be two $F$-stable
$A$-modules. Then $M_1^{F_1}\cong {M_2}^{F_2}$ as $A^F$-modules if
and only if $M_1\cong M_2$ as $A$-modules. In particular, if
$M_1=M_2=M$, then $M^{F_1}\cong M^{F_2}$ as $A^F$-modules.
\end{lem}

\begin{pf} Since $M_i=M_i^{F_i}\otimes_{\bbf_q}k$ and
$A=A^{F}\otimes_{\bbf_q}k$, the lemma follows directly from
Noether-Deuring Theorem (see for example \cite[p.139]{CR}).
\end{pf}

This result shows that, up to isomorphism, it doesn't matter which
Frobenius maps (or $\BF_q$-structures) on $M$ we shall work with
when considering $F$-stable modules.

Let $\modf A$ denote the category whose objects are $F$-stable
modules $M=(M,F_M)$. The morphisms from $(M_1,F_1)$ to $(M_2,F_2)$
are given by homomorphisms $\th\in\hom_A(M_1,M_2)$ such that $\th
\circ F_1=F_2\circ \th$, that is, $A$-module homomorphisms
compatible with their $\BF_q$-structures. Clearly, $\modf A$ is a
subcategory of $\modh A$. It is also easy to see that $\modf A$ is
an abelian $\bbf_q$-category.

The next result allows us to embed a module category defined over
a finite field into a category defined over the algebraic closure
of the finite field.

\begin{thm} \label{CATISO} The abelian category $\modf A$ is equivalent to
the category $\modh A^F$ of finite dimensional $A^F$-modules.
\end{thm}

\begin{pf}
Let $M=(M,F_M)$ be an object in $\modf A$. We define
$\Phi(M):=M^F$, which is an $A^F$-module. Now let
$\thz:(M_1,F_2)\ra (M_2,F_2)$ be a morphism in $\modf A$. Since
$\thz\circ F_1=F_2\circ\thz$, $\thz$ induces a map
$\Phi(\thz):M_1^{F_1}\ra{M_2}^{F_2}$ which is obviously an
$A^F$-module homomorphism. This gives a functor $\Phi:\modf
A\ra\mod A^F$.

Conversely, for each $A^F$-module $X$, we set
$\Psi(X)=X\otimes_{\bbf_q}k$ and define a Frobenius map
$$F_{\Psi(X)}:\Psi(X)\lra \Psi(X);\;
x\otimes\lz\lmto x\otimes\lz^q.$$ By defining
$(a\otimes\lz)(x\otimes\mu)=ax\otimes\lz\mu$ for $a\otimes\lz\in
A^{F}\otimes k=A$ and $x\otimes\mu\in\Psi(X)$ and noting
$F_A(a\otimes \la)=a\otimes \la^q$, we obtain an $A$-module
structure on $\Psi(X)$, which is clearly $F$-stable. Further, for
any morphism $f:X_1\ra X_2$ in $\mod A^F$, the map
$$\Psi(f)=f\otimes 1: \Psi(X_1)\lra \Psi(X_2)$$
is obviously an $A$-module homomorphism satisfying $\Psi(f)\circ
F_1=F_2\circ\Psi(f)$, where $F_i=F_{\Psi(X_i)}$. Hence, we obtain
a functor $\Psi: \mod A^F\ra\modf A$.

From the construction, we see easily from Lemma \ref{IND} that
$$\Psi\Phi=1_{\modf A}\;\;\mbox{and}\;\;\Phi\Psi\cong 1_{\mod A^F},$$
where $1_{\modf A}$ and $1_{\mod A^F}$ denote the identity
functors of $\modf A$ and $\mod A^F$, respectively.
\end{pf}

\begin{cor} There is a one-to-one correspondence between
isoclasses of indecomposable $A^F$-modules and isoclasses of
indecomposable $F$-stable  $A$-modules.
\end{cor}

Let $(M,F_M)$ be an $F$-stable $A$-module. For each submodule $N$
of $M$ (not necessarily an $F_M$-stable subspace), the image
$F_MN$ is an $A$-submodule of $M$.

\begin{prop}\label{soc} Let $(M,F_M)$ and $(N,F_N)$ be two $F$-stable
modules.
\begin{itemize}
\item[(a)] Every submodule $M'$ of $M$ which is also an
$F_M$-stable subspace is an $F$-stable $A$-module. In particular,
both the radical $\rad M$ and socle $\soc M$ of $M$ are $F$-stable
modules.
\item[(b)] As $\BF_q$-spaces, we have
$\Hom_A(M,N)^F\cong \Hom_{A^F}(M^F,N^F)$.
\item[(c)] We have $\BF_q$-algebra isomorphisms
$\End_A(M)^F\cong \End_{A^F}(M^F)$ and
$$\bigl(\End_A(M)/\rad\End_A(M)\bigr)^F\cong
\End_{A^F}(M^F)/\rad\End_{A^F}(M^F).$$ \end{itemize}
\end{prop}

\begin{pf} Since $F_MM'=M'$, the restriction of $F_M$ to $M'$ defines a Frobenius
map on $M'$. So $M'$ is an $F$-stable module as the condition
(\ref{fstab}) is automatically satisfied . If $S$ is a maximal
(resp. simple) submodule of $M$, so is $F_MS$. Thus, both $\rad M$
and $\soc M$ are submodules which are $F_M$-stable (subspaces) and
hence, are $F$-stable modules. (b) is a consequence of the
category isomorphism given in \ref{CATISO}. The first statement in
(c) follows from (b) and Lemma \ref{Lang}. We now prove the last
isomorphism. We first observe that if $B$ is a semisimple algebra
with Frobenius morphism $F_B$, then the fixed point algebra $B^F$
is also semisimple. Since
$$\bigl(\End_A(M)/\rad\End_A(M)\bigr)^F\cong
(\End_A(M))^F/(\rad\End_A(M))^F,$$ it remains to prove that
$$(\rad\End_A(M))^F=\rad(\End_A(M))^F.$$
The inclusion ``$\supseteq$'' follows from the semisimplicity of
the right hand side of the above isomorphism, while the inverse
inclusion ``$\subseteq$'' follows from the fact that
$\bigl(\rad\End_A(M)\bigr)^F$ is nilpotent.
\end{pf}

An $A$-module $M$ is called {\it $F$-periodic}, if there exists an
$F$-stable $A$-module $\tilde M$ such that $M$ is isomorphic to a
direct summand of $\tilde M$ (denoted $M\mid \tilde M$). We shall
see in the next section that for a finite dimensional algebra $A$
with Frobenius morphism $F_A$ every $A$-module is $F$-periodic.
However, in example \ref{nonper}, we shall see that,  for an
infinite dimensional algebra, there are modules which are not
$F$-periodic.

We end this section with the following example which is important
in sections 6--10.

\begin{example}\label{adquiver}
Let $Q=(Q_0,Q_1)$ be a finite quiver without loops, where $Q_0$
resp. $Q_1$ denotes the set of vertices resp. arrows of $Q$ . For
each arrow $\rho$ in $Q_1$, we denote by $h\rho$ and $t\rho$ the
head and the tail of $\rho$, respectively. Let $\sz$ be an
automorphism of $Q$, that is, $\si$ is a permutation on the
vertices of $Q$ and on the arrows of $Q$ such that
$\si(h\rho)=h\si(\rho)$ and $\si(t\rho)=t\si(\rho)$ for any
$\rho\in Q_1$.  We further assume, following \cite[12.1.1]{L93},
that $\si$ is admissible, that is, there are no arrows connecting
vertices in the same orbit of $\sz$ in $Q_0$. We shall call
the pair $(Q,\si)$ an {\it admissible quiver}, or simply an {\it
ad-quiver}.

Let $A:= k Q$ be the path algebra of $Q$ over $k=\ofq$ with
identity $1=\sum_{i\in Q_0}e_i$ where $e_i$ is the idempotent (or
the length 0 path) corresponding to the vertex $i$. Then $\si$
induces a Frobenius morphism (cf. Lemma \ref{frob})
\begin{equation} \label{eq:adquiver}
F_{Q,\si}=F_{Q,\si;q}:A\ra A;\,\,\sum_{s}x_sp_s\lmto
\sum_{s}x_s^q\sz(p_s),
\end{equation} where $\sum_{s}x_sp_s$ is a
$ k$-linear combination of paths $p_s$, and
$\sz(p_s)=\sz(\rho_t)\cdots\sz(\rho_1)$ if
$p_s=\rho_t\cdots\rho_1$ for arrows $\rho_1,\cdots,\rho_t$ in
$Q_1$. We shall investigate the structure of $A^F$ in \S6.
\end{example}


\section{Twisting modules with Frobenius maps}

Let $M$ be an $A$-module and $F_M:M\to M$ a Frobenius map on the
space $M$. Note that $M$ is not necessarily $F$-stable. We define
its {\it (external) Frobenius twist} (with respect to the
Frobenius morphism $F_A$ on $A$) to be the $A$-module $M^{[F_M]}$
such that $M^{[F_M]}=M$ as vector spaces with $F$-twisted action
$$a*m:=F_M\bigl(F_A^{-1}(a)F_M^{-1}(m)\bigr)\text{ for all }a\in A, m\in M.$$
If $\pi:A\to\End_k(M)$ and $\pi^{[F_M]}:A\to\End_k(M^{[F_M]})$
denote the corresponding representations, then
$$\pi^{[F_M]}(a)=F_{(M,M)}(\pi (F_A^{-1}(a)))
=F_M\circ\pi(F_A^{-1}(a))\circ F_M^{-1}\text{ for all }a\in A,$$
where $F_{(M,M)}$ is the induced Frobenius map on $\End_k(M)$ (cf.
\ref{Lang1}).

\begin{lem}  \label{indptF} Up to isomorphism, the Frobenius
twist $M^{[F_M]}$ is independent of the choice of the Frobenius
map $F_M$ on $M$.
\end{lem}
\begin{pf}
If $F_M$ and $F'_M$ are two Frobenius maps on $M$, then the linear
isomorphism $f:=F'_M\circ F_M^{-1}:M\to M$ is clearly an
$A$-module isomorphism from $M^{[F_M]}$ to $M^{[F'_M]}$.
\end{pf}

By the lemma, we shall denote $M^{[F_M]}$ and $\pi^{[F_M]}$ by
$M^{[1]}$ and $\pi^{[1]}$, respectively. Similarly, we define
$M^{[-1]}$ to be the $A$-module given by
$\pi^{[-1]}:A\to\End_k(M)$ where
\begin{equation}\label{inv}\pi^{[-1]}(a) =F_M^{-1}\circ\pi F_A(a)\circ
F_M\text{ for all }a\in A.\end{equation}
Inductively, for each
integer $s> 1$, we define $M^{[s]}=(M^{[s-1]})^{[1]}$ and
$M^{[-s]}=(M^{[-s+1]})^{[-1]}$ with respect to the {\it same}
given $F_M$.

\begin{cor} \label{twist2} The Frobenius twist $(\,\,)^{[1]}$ defines
a category isomorphism from $\modh A$ onto itself.
\end{cor}
\begin{pf} Given two $A$-modules $M$ and $N$ with Frobenius maps $F_M$ and
$F_N$, respectively, and an $A$-module homomorphism $f:M\to N$, it
can be checked by using the corresponding representations that the
linear map $F_{(M,N)}(f)$ defined before Lemma \ref{Lang} is in
fact an $A$-module homomorphism from $M^{[1]}$ to $N^{[1]}$. We
shall denote this morphism by $f^{[1]}$. Thus, we obtain a functor
${(\,\,)}^{[1]}:\modh A\to\modh A$. This functor is clearly
invertible with inverse $(\,\,)^{[-1]}$.
\end{pf}

If $(M,F_M)$ be an $F$-stable $A$-module and $N$ is a submodule of
$M$ (not necessarily an $F_M$-stable space), then the
$A$-submodule $F_MN$ of $M$ is called  the ``{\it internal}''
Frobenius twist of $N$. Note that $F_MN$ is isomorphic to the
(external) Frobenius twist $N^{[1]}$ of $N$ with respect to any
given Frobenius map $F_N$ on $N$. This is deduced from the fact
that the $k$-linear map $\varphi=F_M|_N\circ F_N^{-1} :N^{[1]}\ra
F_MN$ is an $A$-module isomorphism. Recall that $[M]$ denotes the
isoclass of $M$.

\begin{prop}\label{twist} Let $M, M_1$ and $M_2$ be $A$-modules with Frobenius
maps $F_M,F_1$ and $F_2$, respectively.
\begin{itemize}
\item[(a)] $(M,F_M)$ is $F$-stable if and only if, as $A$-modules,
$M^{[1]}=M$.

\item[(b)] $M^{[1]}\cong M$ if and only if there exists $M'\in[M]$ such
that $M'=M$ as a vector space and $(M',F_M)$ is $F$-stable.

\item[(c)] For any given integer $s$,
$M_1\cong M_2$ if and only if $M_1^{[s]}\cong {M_2}^{[s]}$.
\end{itemize}
\end{prop}

\begin{pf}
The statement (a) follows directly from the definition, since
$(M,F_M)$ is $F$-stable if and only if $\pi\circ
F_A=F_{(M,M)}\circ \pi$ which is equivalent to $\pi=\pi^{[F_M]}$ .

We now prove (b). Let $F_M$ be the Frobenius map on $M$ which
defines $M^{[1]}$. Suppose there exists an $f\in GL(M)$ such that
$f\circ\pi^{[1]}(a)=\pi(a)\circ f$ for all $a\in A$. By
Lang-Steinberg's theorem, there exists $g\in GL(M)$ such that
$f=g^{-1}\circ F(g)$, where $F$ is the restriction of $F_{(M,M)}$
on $\End_k(M)$ to $GL(M)$ . Since
$\pi^{[1]}(a)=F_{(M,M)}(\pi(F_A^{-1}(a)))$, it follows that
$F(g)\circ F_{(M,M)}(\pi(F_A^{-1}(a)))\circ F(g^{-1})=g\circ
\pi(a)\circ g^{-1}.$ Putting $a=F_A(b)$, we obtain
$F_{(M,M)}(\pi'(b))=\pi'(F_A(b))$, where $\pi':A\to\End_k(M)$ is
the representation defined by $\pi'(a)=g\circ \pi(a)\circ g^{-1}$
for all $a\in A$. Thus, the module $M'$ (with the same space as
$M$) defined by $\pi'$ is the required one.

To prove (c), by induction, it suffices to prove the cases for
$s=\pm 1$. We only prove the case $s=1$; the proof for the case
$s=-1$ is similar. Let $\pi_1$ and $\pi_2$ be the representations
corresponding to $A$-modules $M_1$ and $M_2$, respectively. Let
$\vphi:M_1\ra M_2$ be a $k$-linear isomorphism. Then
$\psi=F_2\circ\vphi\circ F_1^{-1}$ is also a $k$-linear
isomorphism, where $F_i$ is a Frobenius map on $M_i$. Then,
$\vphi$ is an $A$-module isomorphism if and only if
$\vphi\circ\pi_1(a)=\pi_2(a)\circ\vphi$ for all $a\in A$. Clearly,
the latter is equivalent to
$\psi\circ\pi_1^{[1]}(a)={\pi_2}^{[1]}(a)\circ\psi$ which means
that $\psi:M_1^{[1]}\to M_2^{[1]}$ is an $A$-module isomorphism.
\end{pf}

We now characterize $F$-periodic modules defined at the end of
last section by Frobenius twisting.

\begin{thm} \label{periodic} An $A$-module $M$ is $F$-periodic if
and only if $M^{[r]}\cong M$ for some integer $r$.
\end{thm}

\begin{pf} Suppose $M$ is $F$-periodic. Then so is every direct summand of $M$.
Since $M^{[r]}\cong M$ implies $M^{[rs]}\cong M$ for all $s\ge 1$,
it suffices to prove the case when $M$ is indecomposable.

If $M$ is $F$-periodic and indecomposable, then there is an
$F$-stable indecomposable $A$-module $(N,F_N)$ such that $M\mid
N$. Thus, $F_N^nM\mid N$ for all $n\ge1$. By the
Krull-Remak-Schmidt theorem, there must be a number $r$ such that
$F_N^rM\cong M$, i.e., $M^{[r]}\cong M$.

Conversely, suppose $M^{[r]}\cong M$. By Proposition
\ref{twist}(b), there exists $M'\in[M]$ such that $(M',F_M)$ is an
$F^r$-stable module (with respect to the Frobenius morphisms
$F_A^r$ on $A$ and $F_M^r$ on $M'$), i.e.,  as an $A$-module
${M'}^{[r]}=M'$. Let $\pi:A\to\End_k(M')$ be the corresponding
representation. Then $\pi^{[r]}=\pi$. Set
$$N=M'\oplus {M'}^{[1]}\oplus\cdots\oplus {M'}^{[r-1]}$$
and define a Frobenius map $F_N:N\ra N$ by
\begin{equation}\label{eq:tiF}
F_N(x_0,x_1,\cdots,x_{r-1})
=(F_M(x_{r-1}),F_M(x_0),\cdots,F_M(x_{r-2})).
\end{equation}
Since $\pi^{[r]}=\pi$, it follows that, for any $a\in A$ and
$m=(x_0,x_1,\cdots,x_{r-1})\in N$
$$\aligned
F_N(am)&=F_N(\pi(a)x_0,\pi^{[1]}(a)x_1,\cdots,\pi^{[r-1]}(a)x_{r-1})\cr
&=(\pi^{[r]}(F_A(a))F_M(x_{r-1}),\pi^{[1]}(F_A(a))F_M(x_0),\cdots,
\pi^{[r-1]}(F_A£¨a£©)F_M(x_{r-2}))\cr &=
F_A(a)F_N(m).\endaligned$$ Therefore, $N$ is $F$-stable and,
consequently, $M$ is $F$-periodic.
\end{pf}

\begin{cor} \label{fda} Let $A$ be finite dimensional. Then
every $A$-module is $F$-periodic.
\end{cor}

\begin{pf} Let $M$ be an $A$-module with Frobenius map $F$.
By Lemma \ref{frob}, there are $k$-basis
$\{a_1,a_2,\cdots a_s\}$ of $A$ and $k$-basis
$\{m_1,m_2,\cdots,m_t\}$ of $M$ such that
$$F_A(\sum_{i=1}^s x_ia_i)=\sum_{i=1}^s x_i^q a_i \;\;\text{and}\;\;
F(\sum_{j=1}^ty_jm_j)=\sum_{j=1}^t y_j^q m_j$$ for all
$a=\sum_{i=1}^s x_ia_i\in A$ and all $m=\sum_{j=1}^ty_jm_j$ in
$M$. Write $a_im_j=\sum_{l=1}^t z_{ijl}m_l$ with $z_{ijl}\in k.$
Since $k=\ofq$, there is an integer $n$ such that all $z_{ijl}$'s
lie in $\bbf_{q^n}$. Thus, $F^n(a_im_j)=a_im_j$ for all $i$ and
$j$. Consequently, we have
$$\aligned
F^n(am)&=\sum_{i,j}F^n(x_iy_ja_im_j)
=\sum_{i,j}x_i^{q^n}y_j^{q^n}F^n(a_im_j)\cr &=(\sum_{i=1}^s
x_i^{q^n}a_i)(\sum_{j=1}^t y_j^{q^n}m_j) =F^n_A(a)F^n(m).
\endaligned$$
This proves that $M$ is $F^n$-stable. Now the result follows from
Prop. \ref{twist}(a) and the theorem above.
\end{pf}

Recall that, for a path algebra $A$ of a quiver $Q$, an $A$-module
can be identified as a representation $(V,\phi)$ of $Q$ where
$V=\{V_i\}_{i\in Q_0}$ is a set of finite dimensional vector
spaces $V_i$ and $\phi=\{\phi_\rho\}_{\rho\in Q_1}$ is a set of
linear transformations $\phi_\rho:V_{t\rho}\to V_{h\rho}$.

\begin{example} \label{nonper} Let $Q$ be the quiver with two vertices 1 and 2 and
with infinite many arrows from 1 to 2 indexed by $\rho_{ij}$ for
all $i\geq 1$ and $1\leq j\leq i$. That is, they are
$$\rho_{11}, \rho_{21}, \rho_{22}, \rho_{31},\rho_{32},
\rho_{33},\cdots.$$ Then the path algebra $A=kQ$ is infinite
dimensional, but has the identity $1=e_1+e_2$. Consider the
automorphism $\sigma$ of $Q$ fixing two vertices and cyclicly
permuting each subset $\{\rho_{ij}| 1\leq j\leq i\}$ of arrows for
each fixed $i\geq 1$. Then $\sigma$ induces a Frobenius morphism
$F=F_{Q,\si}$ on $A$ (see \ref{adquiver}). Define a representation
$V$ such that $V_1=V_2=k$ and that $\phi_{\rho_{i1}}$ is the
identity on $k$, but $\phi_{\rho_{ij}}=0$ for all $j\not=1$. Then
$V$ is two-dimensional and $V^{[r]}\not\cong V$ for all $r\ge 0$.
Therefore, $V$ is not $F$-periodic.
\end{example}

\section{$F$-periods and indecomposable $F$-stable  modules}

Let $A$ be a $k$-algebra with a fixed Frobenius morphism $F_A$.
The proof of Theorem \ref{periodic} actually suggests a
construction of indecomposable $F$-stable modules from
$F$-periodic indecomposable ones.

Let $M$ be $F$-periodic with respect to a given Frobenius map
$F_M$, and let $r$ be the minimal integer such that $M^{[r]}\cong
M$. We shall call $r$ the $F$-{\it period} of $M$, denoted
$p(M)=p_F(M)$. Clearly, if $M^{[s]}\cong M$, then $p(M)\mid s$
and, by Lemma \ref{indptF}, $p(M)$ is independent of the choice of
$F_M$.

For $r=p(M)$, by Proposition \ref{twist}(b), there is a new
$A$-module structure $M'$ on the vector space $M$ such that
$M'\cong M$ and $(M',F_M)$ is $F_M^r$-stable, that is,
${M'}^{[r]}=M'$. Thus, by Proposition \ref{twist}(c),
$M',{M'}^{[1]}, \cdots$, ${M'}^{[r-1]}$ are pairwise
non-isomorphic. Let
\begin{equation}\label{Mtil}
{\tilde M}=M'\oplus {M'}^{[1]}\oplus\cdots\oplus {M'}^{[r-1]}
\end{equation} and define a Frobenius map $F_\tM:\tilde M\ra\tilde M$ as in
(\ref{eq:tiF}). Then $({\tilde M},F_\tM)$ is $F$-stable, that is,
$\tilde M\in\Ob(\modf A)$. Thus, $\tM^F:=\tM^{F_\tM}$ is an
$A^F$-module. By Lemma \ref{IND} we infer that, up to isomorphism,
${\tilde M}^F$ is independent of the choice of $M'$.

The following result generalizes Kac's result \cite[Lemma 3.4]{K1}
for the path algebras of quivers.

\begin{thm} \label{INDF} Maintain the notation above. Let $M$ be an $F$-periodic
indecomposable $A$-module with $F$-period $r$. Then $(\tilde M,
F_\tM)$ is indecomposable in $\modf A$ and
$$\End_{A^F}({\tilde M}^F)/\rad(\End_{A^F}({\tilde M}^F))
\cong\bbf_{q^r}.$$ Moreover, every indecomposable $A^F$-module is
isomorphic to a module of the form ${\tilde M}^F$ for
some $F$-periodic indecomposable module $M$.
\end{thm}

\begin{pf} The Frobenius map $F_\tM:{\tilde M}\ra {\tilde M}$ induces
a Frobenius map ${\tilde F}=F_{(\tM,\tM)}$ on $\End_A({\tilde M})$.
By Lemma \ref{soc}, we have an $\bbf_q$-algebra isomorphism
$$\End_{A^F}({\tilde M}^F)/\rad(\End_{A^F}({\tilde M}^F))
\cong \bigl(\End_A(\tilde M)/\rad(\End_A(\tilde M))\bigr)^{\tilde
F}.$$ For each $f\in \End_A(\tilde M)$, we write $f$ in matrix
form as
$$f=(f_{ji})_{r\times r}:\tilde M=\bigoplus_{i=0}^{r-1} {M'}^{[i]}
\lra \bigoplus_{i=0}^{r-1}{M'}^{[i]}=\tilde M,$$ where
$f_{ji}:{M'}^{[i]}\ra {M'}^{[j]}$. Then ${\tilde F}(f)=(g_{ji})$
with $g_{ji}=f_{j-1,i-1}^{[1]}$, where the indices are integers
modulo $r$ . In particular, if ${\tilde F}(f)=f$, i.e.,
$f\in\End_A(\tilde M)^{\tilde F}$, then $f_{ji}=f_{j-1,i-1}^{[1]}$
for all $0\leq i,j\leq r-1$. Since $M', {M'}^{[1]},\cdots,
{M'}^{[r-1]}$ are pairwise non-isomorphic indecomposable
$A$-modules, we have an algebra isomorphism
$$\aligned
\End_A(\tilde M)/\rad(\End_A(\tilde M)) &\lra
\prod_{i=0}^{r-1}\End_A({M'}^{[i]})/\rad(\End_A({M'}^{[i]}))\cr
{\bar f}=({\bar f}_{ji}) &\lmto ({\bar f}_{00},\bar
f_{11}\cdots,{\bar f}_{r-1,r-1}).\cr
\endaligned$$
Since
$\End_A({M'}^{[i]})/\rad(\End_A({M'}^{[i]}))\cong k$ for each
$0\leq i\leq r-1$, we obtain
$$\End_A(\tilde M)/\rad(\End_A(\tilde M))
\cong \underbrace{k\times\cdots \times k}_{r}=:U.$$
The Frobenius map $\tilde F$ on the left hand side induces a Frobenius map
$\tilde F$ on $U$ given by
$$\tilde F(x_0,x_1,\cdots,x_{r-1})=(x_{r-1}^q,x_0^q,\cdots,x_{r-2}^q).$$
Hence,
$$\End_{A^F}({\tilde M}^F)/\rad(\End_{A^F}({\tilde M}^F))
\cong U^{\tilde F}\cong \bbf_{q^r}.$$

Conversely, let $X$ be an indecomposable $A^F$-module such that
$$\End_{A^F}(X)/\rad(\End_{A^F}(X))\cong \bbf_{q^r}.$$
Then $\bX:=X\otimes_{\bbf_q}k$ is an $F$-stable $A$-module
with the Frobenius map $F:=F_\bX$ defined by
$$F(x\otimes\lz)=x\otimes \lz^q\;\;
\text{for}\;\, x\otimes\lz\in X\otimes_{\bbf_q}k.$$
Moreover, we have a decomposition
$$\bX=M_1\oplus\cdots\oplus M_r,$$
where $M_1,\cdots,M_r$ are pairwise non-isomorphic indecomposable
$A$-modules. The $F_A$-stability of $\bX$ implies that
$$F M_1\oplus\cdots\oplus F M_r=F{\bX}
=\bX=M_1\oplus\cdots\oplus M_r.$$ Thus, the set
$\ooz:=\{M_1,\cdots,M_r\}$ is $F$-stable (up to isomorphism). We
claim that $\ooz$ contains only one $F$-orbit. Let
$\ooz_1,\cdots,\ooz_m$ be the orbits of $\ooz$. For each $1\leq
j\leq m$, let
$${\hat M_j}=\bigoplus_{N\in \ooz_j}N.$$
Then
$${\hat M}_j=N_{j}\oplus F N_{j}\oplus\cdots\oplus
F^{r_j-1} N_{j},$$ where $N_{j}\in \ooz_j$ and $r_j=|\ooz_j|$.
Note that $F^{r_j} N_{j}\cong N_{j}$. Let $F_j$ be a Frobenius map
on $N_{j}$. Then $FN_{j}$ is isomorphic to the Frobenius twist
$N_{j}^{[1]}$ of $N_{j}$ with respect to $F_j$. Hence, we have
$N_{j}^{[r_j]}\cong N_{j}$ and
$${\hat M}_j\cong N_{j}\oplus N_{j}^{[1]}\oplus\cdots\oplus
N_{j}^{[r_j-1]}\cong {\hat M}_j^{[1]}.$$ Then, by Prop.
\ref{twist}(b), we obtain an $F$-stable module $({\tilde M}_j,
\tilde F_{j})$ satisfying ${\tilde M}_j\cong {\hat M}_j$. The
Frobenius maps $\tilde F_{j}$ induces a Frobenius map $F'$ on
$\oplus_{j=1}^m {\tilde M}_j$ defined by
$$F'(v_1,\cdots,v_m)=(\tilde F_{1}(v_1),\cdots, \tilde F_{m}(v_m)).$$
By Lemma \ref{IND}, the isomorphism $\bX\cong \oplus_{j=1}^m {\tilde M}_j$
implies that
$$X={\bX}^F\cong (\oplus_{j=1}^m {\tilde M}_j)^{F'}
=\oplus_{j=1}^m {\tilde M}_j^{\tilde F_{j}}.$$ Since $X$ is
indecomposable, we must have $m=1$ and the required isomorphism
$$\bX= {\hat M}_1\cong N_{1}\oplus N_{1}^{[1]}\oplus\cdots\oplus
N_{1}^{[r-1]}.$$
\end{pf}

From the proof, we obtain the following correspondence.

\begin{cor} If $A$ is finite dimensional, then there is a
one-to-one correspondence between the isoclasses of indecomposable
$A^F$-modules and the $F$-orbits of the isoclasses of
indecomposable $A$-modules.
\end{cor}

 A finite dimensional algebra $B$ over a field is called
representation-finite if there are only finitely many isoclasses
of (finite dimensional) indecomposable $B$-modules. Theorem
\ref{INDF} and Cor. \ref{fda} imply immediately the following.

\begin{cor} Let $A$ be finite dimensional and $F_A$ a Frobenius morphism on
$A$. Then, $A$ is representation-finite if and only if so is
$A^F$.
\end{cor}

An $A^F$-module $X$ is called absolutely indecomposable if
$X\otimes_{\bbf_q}k$ is an indecomposable $A$-module.

\begin{cor} An $A^F$-module $X$ is absolutely indecomposable if and only if
there is an $F$-stable indecomposable $A$-module $M=(M,F)$ such
that $X\cong M^F$.
\end{cor}

\section{Finite dimensional hereditary algebras}

The first application of our theory is to show that every finite
dimensional  hereditary (basic) algebra over a finite field is
isomorphic to the $F$-fixed point algebra of the path algebra of a
finite ad-quiver (see Example \ref{adquiver}). Thus, the
representation theory of a finite dimensional hereditary algebra
is  completely determined by the counterpart of the corresponding
ad-quivers.

We first recall the notion of modulated quivers (cf.
\cite{DR1} and \cite[4.1.9]{Be}).

\begin{defn}\label{vq}  A valued graph  is a graph without loops
together with a positive
integer $d_x$ for each vertex $x$ and a pair of positive integers
$({}_xc_y^\ga, {}_yc_x^\ga)$ for each edge $x\,\overset
\ga{\text{---}}\,y$ satisfying ${}_xc_y^\ga d_y={}_yc_x^\ga d_x$.
A valued graph together with an orientation is called a {\it
valued quiver}, and a valued quiver is called {\it simple} if it
has no parallel arrows. (Thus, opposite arrows between two
vertices are allowed in a simple valued quiver.) A {\it
modulation}\footnote{If we follow the definition given in
\cite[4.1.9]{Be}, then $M_\rho=\begin{cases}{}_xM_y^\ga,&\text{ if
}\rho=y\overset \ga\longrightarrow x\cr {}_yM_x^\ga,&\text{ if
}\rho=x\overset \ga\longrightarrow y.\cr
\end{cases}$}
$\BM$ of a valued quiver consists of an assignment of
a division ring $D_x$ to each vertex $x$, and a
$D_y$-$D_x$-bimodule $M_\rho$ to each arrow $\rho=x\longrightarrow
y$ satisfying
\begin{itemize}
\item[(1)] $\hom_{D_y}(M_\rho,D_y)\cong
\hom_{D_x}(M_\rho,D_x)$,
\item[(2)] $\dim_{D_y}(M_\rho)={}_xc_y^\ga$, $\dim_{D_x}(M_\rho)={}_yc_x^\ga$.
\end{itemize}
Finally, a {\it modulated quiver} consists of a valued quiver
$\Ga$ and a modulation $\BM$.

A modulated quiver is simple if its underlying quiver is simple.
\end{defn}

Let ${\cal Q}=(\Ga,\BM)$ be a modulated quiver with $\Ga_0$ (resp.
$\Ga_1$) the set of vertices (resp. arrows) of $\Ga$ and
$\BM=(\{D_x\}_{x\in\Ga_0}, \{M_\rho\}_{\rho\in\Ga_1})$. Let
$R=\oplus_{x\in\Ga_0}D_x$ and  $M=\oplus_{\rho\in\Ga_1}M^\rho$.
Then $M$ is a natural $R$-$R$-bimodules. The $R$-algebra
$$T({\cal Q}):=\bigoplus_{n\ge 0}M^{\otimes n}\text{ where } M^{\otimes0}=R,
M^{\otimes n}=\underbrace{M\otimes_R\cdots\otimes_RM}_n$$ is
called the {\it path} (or {\it tensor}) algebra of $\cal Q$. Thus,
a tensor $x_n\otimes\cdots\otimes x_1$ with  $x_i\in M_{\rho_i}$
is non-zero implies that $\rho_n\cdots \rho_1$ is a path in $\Ga$.

For a modulated quiver ${\cal Q}=(\Ga,\BM)$, let $\bar{\cal
Q}=(\bar\Ga,\bar\BM)$ be the associated simple modulated quiver
obtained by summing the valuations and bimodules over parallel
arrows. More precisely, $\bar\Ga$ is a simple valued quiver
defined by setting $\bar\Ga_0=\Ga$, $\bar\Ga_1=\{\bar\rho:x\to
y\}$ where ${\bar\rho}:=\{\text{ all }\rho:x\to y\text{ in
}\Ga\}$, and setting the valuation  for the arrow $\bar\rho:x\to
y$ to be $(c_{\bar\rho},c'_{\bar\rho})$ where
$c_{\bar\rho}=\sum_{\rho\in\bar\rho}{}_xc_y^\rho$ and
$c_{\bar\rho}'=\sum_{\rho\in\bar\rho}{}_yc_x^\rho$. The modulation
$\bar{\BM}=(\{D_x\}_{x\in\bar\Ga_0},\{\bar
M_{\bar\rho}\}_{\bar\rho\in\bar\Ga_1})$ is defined by setting
$\bar M_{\bar\rho}=\oplus_{\rho\in\bar\rho}M_\rho$.

\begin{defn} \label{mqi}
Let ${\cal Q}=(\Ga,\BM)$ and ${\cal Q}'=(\Ga',\BM')$ be two
modulated quiver. We say that ${\cal Q}\cong {\cal Q}'$ if there
exists a quiver isomorphism
$\tau:\bar\Ga\overset\sim\longrightarrow\bar\Ga'$ such that (1)
$D_x\cong D'_{\tau(x)}$ as division rings, and (2)
$M_{\bar\rho}\cong M'_{\tau(\bar\rho)}$ as bimodules via (1).
\end{defn}

Clearly, if ${\cal Q}\cong {\cal Q}'$ then we have algebra
isomorphism $T({\cal Q})\cong T({\cal Q}')$.

We now construct a modulated quiver from an ad-quiver.  Given a
finite ad-quiver $(Q,\si)$, let $I=\Ga_0$ and $\Ga_1$ denote the
set of $\si$-orbits in $Q_0$ and $Q_1$, respectively. Thus, we
obtain a new quiver $\Ga=(\Ga_0,\Ga_1)$. For each arrow
$\brho:\bi\longrightarrow \bj$ in $\Ga$, define
\begin{equation}\label{GCM}
\ez_{\brho}=\#\{\mbox{arrows in $\brho$}\},
\;d_\brho=\ez_\brho/\ez_\bj,\;\;
\text{and}\;\;d'_{\brho}=\ez_\brho/\ez_\bi,
\end{equation} where
$\ez_\bk=\#\{\text{vertices in $\si$-orbit } \bk\}$ for $\bk\in
I$. The quiver $\Ga$ together with the valuation
$\{\ez_\bd\}_{\bd\in\Ga_0}$,
$\{(d_{\brho},d'_{\brho})\}_{\brho\in\Ga_1}$ defines a valued
quiver
 $\ggz=\ggz(Q,\sz)$.
Clearly, each valued quiver can be obtained in this way from an
ad-quiver.

Using the Frobenius morphism $F=F_{Q,\si}$ on $A$ defined in
(\ref{eq:adquiver}), we can attach naturally  to $\Ga$ an
$\bbf_q$-modulation  to obtain a modulated quiver (i.e., an
$\BF_q$-species) as follows:
for each vertex $\bi\in I$ and each arrow $\brho$ in $\ggz$, we
fix $i_0\in\bi$, $\rho_0\in \brho$, and consider the $F_A$-stable
subspaces of $A$
$$A_\bi=\bigoplus_{i\in\bi} k e_i=\bigoplus_{s=0}^{\ez_\bi-1} k e_{\sz^s(i_0)}
\;\mbox{\ and\ }\; A_\brho=\bigoplus_{\rho\in\brho} k
\rho=\bigoplus_{t=0}^{\ez_\brho-1} k\sz^t(\rho_0),$$ where $e_i$
denotes the idempotent corresponding to the vertex $i$. Then
\begin{equation}\label{eq:AiF}
A_\bi^F=\{\sum_{s=0}^{\ez_\bi-1}x^{q^s}e_{\sz^s(i_0)}\mid x\in k,
x^{q^{\ez_\bi}}=x\} \text{ and }
A_\brho^F=\{\sum_{t=0}^{\ez_\brho-1}x^{q^t}\sz^t(\rho_0) \mid x\in
k, x^{q^{\ez_\brho}}=x\}. \end{equation}  Further, the algebra
structure of $A$ induces an $A_\bj^F$-$A_\bi^F$-bimodule structure
on $A_\brho^F$ where $\brho:\bi\longrightarrow \bj$. Thus, we
obtain an $\bbf_q$-modulation
$\BM=\BM(Q,\si):=(\{A_\bi^F\}_\bi,\{A_\brho^F\}_\brho)$ over the
valued quiver $\Ga$. We shall denote the $\BF_q$-modulated quiver
${\cal Q}=(\Ga,\BM)$ defined above by
\begin{equation}\label{mqs}
\fM_{Q,\sz}=\fM_{Q,\sz;q}=(\Ga,\BM).
\end{equation}

Let $T(\fM_{Q,\sz})$ be the tensor algebra of the modulated quiver
$\fM_{Q,\si}$. Thus, by definition, $T(\fM_{Q,\sz})=\oplus_{n\ge
0}M^{\otimes n}$, where $M=\oplus_{\brho\in\Ga_1} A_\brho^F$ is
viewed as an $R$-$R$-bimodule with $R=\oplus_{\bi\in I}A_\bi^F$
and $\otimes=\otimes_R$. If, for each $\sz$-orbit $\bf p$ of a
path $\rho_n\cdots\rho_2\rho_1$ in $Q$, we set $A_{\bf
p}=\oplus_{p\in \bf p} kp$. Then
$$A_{\bf p}^F\cong A_{\brho_n}^F\otimes_{\BF_{n-1}}
\cdots\otimes_{\BF_2} A_{\brho_2}^F \otimes_{\BF_1}
A_{\brho_1}^F,$$ where $\brho_t$ is the $\si$-orbit of $\rho_t$
and $\BF_t=A_{h\brho_t}^F$. Since $A^F=\oplus_{\bf p}A_{\bf p}^F$,
it follows that the fixed point algebra $A^F$ is isomorphic to the
tensor algebra $T(\fM_{Q,\sz})$. Thus, $A^F$-modules can be
identified with representations of the modulated quiver
$\fM_{Q,\sz}$ (see \cite{DR1}). The above observation together
with Theorem \ref{CATISO} implies the following.

\begin{prop} \label{sigr} Let $(Q,\si)$ be an ad-quiver with path algebra
$A=kQ$ and induced Frobenius morphism $F=F_{Q,\si}$. Let
$\fM_{Q,\si}$ be the associated $\BF_q$-modulated quiver defined
as above.
\begin{itemize}
\item[(a)] We have an algebra isomorphism $A^F\cong T(\fM_{Q,\si})$.
Hence the categories $\modf A$ and $\mod T(\fM_{Q,\sz})$ are
equivalent.
\item [(b)] If $Q$ has no oriented cycles, then the
fixed-point algebra $A^F$ is a finite dimensional hereditary basic
algebra.
\end{itemize}
\end{prop}

\begin{cor} Maintain the notation above and let $r\geq 1$ be an integer.
Then, the ad-quiver $(Q,\si^r)$ defines an $\bbf_{q^r}$-modulated
quiver $\fM_{Q,\sz^r;q^r}$ whose tensor algebra is isomorphic to
the $\bbf_{q^r}$-algebra $A^F\otimes_{\BF_q}\BF_{q^r}$.
\end{cor}

\begin{pf} Clearly, $F_{Q,\si^r;q^r}=F_{Q,\si;q}^r=F^r$. By Lemma
\ref{frob}, we have $A=A^F\otimes_{\bbf_q} \overline{\bbf}_q$ and
$F(a\otimes \lz)=a\otimes \lz^q$. Then
$A^{F^r}=A^F\otimes_{\bbf_q}{\bbf_{q^r}}$. Now, the isomorphism
follows from Prop. \ref{sigr}(a).
\end{pf}

Note that  for $r>1$ the modulated quiver $\fM_{Q,\si^r;q^r}$ is
different from the modulated quivers $\fM_{Q,\si^r;q}$ and
$\fM_{Q,\si;q^r}$. The former has the same underlying valued
quiver as $\fM_{Q,\si^r;q^r}$, but different base field, while the
latter has the same base field but different underlying valued
quiver (if $\si\neq1$).

Our next result shows that the converse of Prop. \ref{sigr} (b) is
also true.

\begin{thm}\label{ALGISO} Let $B$ be a finite dimensional hereditary basic
algebra over $\BF_q$. Then there is an ad-quiver $(Q,\si)$ such that $B$ is
isomorphic to $(kQ)^{F_{Q,\si}}$.
\end{thm}
\def\gr{{\text{\rm gr\,}}}
\begin{pf} Let $A=B\otimes k$ and define $F:A\to A$ by $F(b\otimes
\la)=b\otimes \la^q$. Clearly, $F$ is a Frobenius morphism on the
$k$-algebra $A$ and $A^F= B$. Since $B$ is a finite dimensional
hereditary basic algebra,  it follows that so is $A$ and $B$ is
isomorphic to the tensor (or path) algebra of the associated
modulated Ext-quiver $Q_B$ (see \cite[p.104]{Be}). In particular,
we have algebra isomorphism
$$\ph:B\overset\sim\lra \gr B\text{\quad where\quad} \gr B=\oplus_{i\ge0}
\rad^iB/\rad^{i+1}B.$$ Since $\rad A=\rad(B\otimes k)= (\rad
B)\otimes k$ (see, e.g., \cite[p.146]{CR}), we have by induction
$\rad^iA=(\rad^iB)\otimes k$ for all $i\ge 1$. Thus, the
isomorphism $\ph$ induces a $k$-algebra isomorphism
$$\tilde\ph:A=B\otimes k\overset\sim\longrightarrow\gr A=\gr B\otimes k.$$
Let $\bar F$ denote the Frobenius morphism on $\gr A$ induced from
$F$. Clearly, $\bar F$ stabilizes each direct summand
$\rad^iA/\rad^{i+1}A$ of $\gr A$ and $B=A^F\cong (\gr A)^{\bar
F}.$

We now prove that $(\gr A,\bar F)$ defines an ad-quiver $(Q,\si)$
such that $(\gr A)^{\bar F}\cong (kQ)^{F_{Q,\si}}$. Suppose
$$A/\rad A=k\bar e_1\oplus\cdots\oplus k\bar e_n$$
where $e_i$ are primitive orthogonal idempotents of $A$ with
$1=e_1+\cdots+e_n$. Since the $F(e_i)$'s form a complete set of
primitive orthogonal idempotents of $A$, it follows that there is
a permutation $\si$ and an invertible element $u\in A$ such that
$F(e_i)=ue_{\si i}u^{-1}$. Thus, $\bar F(\bar e_i)=\bar e_{\si i}$
for all $i$ .

Let $Q_0=\{1,2,\cdots, n\}$ and $I$ the set of $\si$-orbits.
Putting $f_\bi=\sum_{i\in\bi}e_i$ for each $\bi\in I$, we have
$$\rad A/\rad^2 A=\bigoplus_{\bi,\bj\in I}f_\bi(\rad A/\rad^2
A)f_\bj,$$ where $$f_\bi(\rad A/\rad^2 A)f_\bj=\bar f_\bi(\rad
A/\rad^2 A)\bar f_\bj=\bigoplus_{i\in\bi,j\in\bj} \bar e_i(\rad
A/\rad^2 A)\bar e_j.$$ Clearly, $\bar F$ stabilizes each summand
$f_\bi(\rad A/\rad^2 A)f_\bj$ and $\bar F(V_{ij})=V_{\si i,\si j}$
where
$$V_{ij}=\bar e_i(\rad A/\rad^2 A)\bar e_j.$$ Fix $i,j$ and
let $s$ be the smallest integer such that $\si^si=i$ and
$\si^sj=j$. Since $\bar F^s$ stabilizes $V_{ij}$,
we can choose a $k$-basis $v_1,\cdots,v_{t_{ij}}$ of $V_{ij}$ such
that $\bar F^s(v_a)=v_a$ for all $1\le a\le t_{ij}$. Thus, we
obtain a $k$-basis
$$\{v_1,\cdots,v_t,\bar F(v_1),\cdots,
\bar F(v_t),\cdots,\bar F^{s-1}(v_1),\cdots, \bar F^{s-1}(v_t)\}$$
for the $\bar F$-stable space
$\oplus_{a=0}^{s-1}V_{\si^ai,\si^aj}$. Clearly, such a basis can
be constructed for every $\si$-orbit of the set $\{(i,j)\mid
i\in\bi,j\in\bj\}$. Thus, we obtain a basis $\{v^{ij}_a\}_{i,j,
a}$ for $f_\bi(\rad A/\rad^2 A)f_\bj$, and hence for $\rad
A/\rad^2 A$, which is $\bar F$-stable, i.e., $\{\bar
F(v^{ij}_a)\}_{i,j,a}=\{v^{ij}_a\}_{i,j,a}$. If $Q_1$ is the set
of arrows $\rho^{ij}_a$ indexed by the basis elements $v^{ij}_a$,
then $\bar F$ induces a permutation $\si$ on $Q_1$.

So we have obtained an ad-quiver $(Q,\si)$. The standard
$k$-algebra isomorphism $\psi: kQ\to \gr A$ sending $i$ to $\bar
e_i$ and $\rho^{ij}_a$ to $v^{ij}_a$ is compatible with the
Frobenius morphisms $\bar F$ and $F_{Q,\si}$, that is, $\psi\circ
F_{Q,\si}=\bar F\circ\psi$. Consequently, we obtain
$$B\cong(\gr A)^{\bar F}\cong(kQ)^{F_{Q,\si}}.$$
\end{pf}

\begin{cor} Every finite dimensional hereditary algebra over a finite
field is Morita equivalent to the $F$-fixed point algebra of the path
algebra of an ad-quiver.
\end{cor}

\begin{rems} \label{Hua-Hub} (a) If we identify $A_\bi^F$ and $A_\brho^F$
defined in (\ref{eq:AiF}) with $\fqi$ and $\fqaz$ via
$$\sum_{i\in\bi}x_ie_i\lmto x_{i_0}\;\;\text{and}\;\;
\sum_{\rho\in\brho}y_\rho\rho\lmto y_{\rho_0},$$ respectively,
where $i_0\in\bi$ and $j_0\in\bj$ are fixed,
 then, for an arrow $\brho:\bi\ra\bj$, the
induced $\fqj$-$\fqi$-bimodule structure on $A_\brho=\fqaz$ is not
necessarily the natural bimodule $_{\fqj}({\fqaz})_{\fqi}$ induced
by the subfield structure. This is because, for $r,s\geq 1$, the
$\bbf_q$-algebra isomorphism
$$\bbf_{q^r}\otimes_{\bbf_q} \bbf_{q^s}\cong
\underbrace{\bbf_{q^m}\times\cdots\times \bbf_{q^m}}_{d},$$ means
that there are exactly $d$ non-isomorphic simple
$\bbf_{q^r}$-$\bbf_{q^s}$-bimodules. Here $d$ and $m$ denote the
greatest common divisor and the least common multiple of $r,s$,
respectively. More precisely, for $0\leq r\leq \ez_\bi-1$ and
$0\leq s\leq\ez_\bj-1$ with $\rho_0:\sz^r(i_0)\ra\sz^s(j_0)$, the
induced ${\fqj}$-${\fqi}$-bimodule structure on $A_\brho=\fqaz$ is
given by
$$y\cdot z\cdot x=y^{q^s}zx^{q^r}\;\;\text{for}\;\;
x\in\fqi, y\in\fqj, z\in\fqaz,$$ while the natural bimodule
structure corresponds only to the case when $r=s=0$. These two
bimodule structures on $\fqaz$ are not necessarily isomorphic.

(b) The $\bbf_q$-modulated quiver studied in \cite{Hua2},
\cite{Hub1}, and \cite{Hub2} involve only the natural bimodules
$_{\bbf_{q^r}}({\bbf_{q^n}})_{\bbf_{q^s}}$, where $r,s,n\geq 1$,
$r\mid n$ and $s\mid n$. These will be called {\it natural
modulated quivers} below. Note that the natural
$\bbf_{q^r}$-$\bbf_{q^s}$-bimodule
$_{\bbf_{q^r}}({\bbf_{q^n}})_{\bbf_{q^s}}$ is isomorphic to the
direct sum of $\frac{n}{m}$ copies of the natural simple bimodule
$_{\bbf_{q^r}}({\bbf_{q^m}})_{\bbf_{q^s}}$. Though
non-isomorphic bimodules involved in two modulated quivers could
result in isomorphic tensor algebras, the following example shows
that not every finite dimensional hereditary basic algebra over
$\bbf_q$ arises from a natural modulated quiver.
\end{rems}

\begin{example} Let $\ggz$ denote the valued quiver
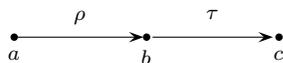
\begin{center}
\begin{pspicture}(0,-0.2)(6,1.2)
\psset{xunit=.8cm,yunit=.7cm,linewidth=.4pt} \psdot*(0.6,0.6)
\psdot*(2.8,0.6) \psdot*(5,0.6) \uput[d](0.6,0.6){$_a$}
\uput[d](2.8,0.6){$_b$} \uput[d](5,0.6){$_c$}
\uput[u](1.7,0.6){$_\rho$} \uput[u](3.9,0.7){$_\tau$}
\psline[arrowsize=4pt]{->}(0.65,0.6)(2.75,0.6)
\psline[arrowsize=4pt]{->}(2.9,0.6)(4.9,0.6)
\end{pspicture}
\end{center}
with $(\ez_a,\ez_b,\ez_c)=(1,2,2)$, $(d_\rho,d'_\rho)=(1,2)$, and
$(d_\tau,d'_\tau)=(2,2)$. Then the pair of natural bimodules
$_{\bbf_q}({\bbf_{q^2}})_{\bbf_{q^2}}$ and
$_{\bbf_{q^2}}({\bbf_{q^4}})_{\bbf_{q^2}}$ together with the
valued quiver $\ggz$ defines the natural $\bbf_q$-modulated quiver
$\fM$. Let $M$ denote the natural bimodule
$_{\bbf_q^2}({\bbf_{q^2}})_{\bbf_{q^2}}$ and $M'=\bbf_{q^2}$ (as
$\bbf_q$-vector spaces) denote the
$\bbf_{q^2}$-$\bbf_{q^2}$-bimodule given by $x\cdot y\cdot
z=xyz^q$. It is easy to see that $M$ and $M'$ are not isomorphic
as $\bbf_{q^2}$-$\bbf_{q^2}$-bimodules. The pair of bimodules
$_{\bbf_q}({\bbf_{q^2}})_{\bbf_{q^2}}$ and $M\oplus M'$ also
defines a new $\bbf_q$-modulated quiver $\fM'$ whose tensor
algebra $T(\fM')$ is not isomorphic to the tensor algebra $T(\fM)$
of $\fM$, since $T(\fM)\otimes_{\bbf q} k$ and
$T(\fM')\otimes_{\bbf q} k $ are respectively isomorphic to the
path algebras of the following quivers
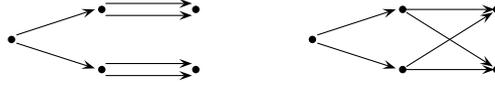
\begin{center}
\begin{pspicture}(-2,0)(10,1.5)
\psset{xunit=.5cm,yunit=.4cm,linewidth=.4pt}
\psdot*(0.6,1.5)
\psdot*(3,0.5) \psdot*(3,2.5) \psdot*(5.5,0.5) \psdot*(5.5,2.5)
\psline[arrowsize=3.5pt]{->}(0.73,1.65)(2.82,2.5)
\psline[arrowsize=3.5pt]{->}(0.73, 1.38)(2.82,0.54)
\psline[arrowsize=3.5pt]{->}(3.1,2.7)(5.4,2.7)
\psline[arrowsize=3.5pt]{->}(3.1,2.3)(5.4,2.3)
\psline[arrowsize=3.5pt]{->}(3.1,0.7)(5.4,0.7)
\psline[arrowsize=3.5pt]{->}(3.1,0.3)(5.4,0.3)

\psdot*(8.6,1.5)
\psdot*(11,0.5) \psdot*(11,2.5) \psdot*(13.5,0.5) \psdot*(13.5,2.5)
\psline[arrowsize=3.5pt]{->}(8.73,1.65)(10.82,2.5)
\psline[arrowsize=3.5pt]{->}(8.73, 1.38)(10.82,0.54)
\psline[arrowsize=3.5pt]{->}(11.05,2.5)(13.45,2.5)
\psline[arrowsize=3.5pt]{->}(11.05,0.5)(13.45,0.5)
\psline[arrowsize=3.5pt]{->}(11.1,0.6)(13.4,2.4)
\psline[arrowsize=3.5pt]{->}(11.1,2.4)(13.4,0.6)
\end{pspicture}
\end{center}
and are obviously not isomorphic.
\end{example}

\section{Almost split sequences}

The Auslander-Reiten theory is one of the fundamental tools in the
study of representations of algebras (see, e.g., \cite{ARS}). This
and next sections are devoted to establishing a relation between
the Auslander-Reiten theories of $A$ and its fixed-point algebra
$A^F$.

We briefly review the general theory. Let $A$ be a finite
dimensional algebra over an {\it arbitrary} field $k$. A morphism
$\vphi:L\ra M$ in $\mod A$ is called {\it minimal right almost
split} if
\begin{itemize}
\item[(a)] $\vphi$ is not a split epimorphism,

\item[(b)] any morphism $X\ra M$ which is not a split epimorphism factors
           through $\vphi$,

\item[(c)] any morphism $f:L\ra L$ satisfying $\vphi=\vphi f$ is an isomorphism.
\end{itemize}
It is easy to see that, if $\vphi:L\ra M$ is a minimal right
almost split morphism, then $M$ is indecomposable.

A {\it minimal left almost split} morphism is defined dually.

A short exact sequence $0\ra
N\skrel{\psi}{\ra}L\skrel{\vphi}{\ra}M\ra 0$ is called an {\it
almost split sequence} if $\vphi$ is minimal right almost split,
or equivalently, if $\psi$ is minimal left almost split.

\begin{thm} \label{AR} (Auslander-Reiten) Let $M$ be an indecomposable $A$-module.
\begin{itemize}
\item[(a)] There exists a unique, up to isomorphism, minimal right almost
split morphism $\vphi:L\ra M$. If, moreover, $M$ is not
projective, then the sequence
$$0\lra\Ker\vphi\lra L\skrel{\vphi}{\lra}M\lra 0$$
is an almost split sequence.

\item[(b)] There exists a unique, up to isomorphism, minimal left almost
split morphism $\psi:M\ra L$. If, moreover, $M$ is not injective,
then the sequence
$$0\lra M\skrel{\psi}{\lra} L\lra \Coker\psi\lra 0$$
is an almost split sequence.
\end{itemize}
\end{thm}

Given an almost split sequence $0\ra
N\skrel{\psi}{\ra}L\skrel{\vphi}{\ra}M\ra 0$, the module $N$
(resp. $M$) is uniquely (up to isomorphism) determined by $M$
(resp. $N$). We write $\tau_A M=N$ or $N=\tau_A^{-1}M$. The
$\tau_A=:\tau$ is called the {\it Auslander-Reiten translation} of
$A$, which indeed admits the following description using the
transpose and the dual (see \cite[Chap. IV, V]{ARS}. Let $\cp$ and
$\ci$ be the ideals of $\mod A$ (in the sense of \cite[p.16]{GR})
defined respectively by
$$\cp(X,Y)=\{f:X\ra Y|\text{$f$ factors through a projective $A$-module}\}$$
and
$$\ci(X,Y)=\{f:X\ra Y|\text{$f$ factors through an injective $A$-module}\},$$
where $X,Y\in\mod A$. The factor categories of $\mod A$ by $\cp$
and $\ci$ are denoted by $\undermod A$ and $\overmod A$,
respectively. Let $M\in\mod A$ and let
$$P_1\skrel{f}{\lra} P_0\lra M\lra 0$$
be a minimal projective presentation of $M$. Applying
$\hom_A(-,A)$, we get a morphism
$$f^\ast:=\hom_A(f,A):\hom_A(P_0,A)\ra \hom_A(P_1,A)$$
in $\mod A^{\op}$. Thus, we obtain an $A^{\op}$-module $\Coker
f^\ast$ which is called the {\it transpose} of $M$ and denoted by
$\Tr M$. In general, $\Tr$ does not give rise to a functor from
$\mod A$ to $\mod A^{\op}$. However, $\Tr$ induces a functor
$\Tr:\undermod A\ra\undermod A^{\op}$ such that the compositions
$\undermod A\skrel{\Tr}{\ra}\undermod A^{\op}\skrel{\Tr}{\ra}
\undermod A$ and $\undermod A^{\op}\skrel{\Tr}{\ra}\undermod
A\skrel{\Tr}{\ra} \undermod A^{\op}$ are isomorphic to the
identity on $\undermod A$ and $\undermod A^{\op}$, respectively.
On the other hand, the duality $D=\hom_k(-,k):\mod A\ra\mod
A^{\op}$ induces a duality $D:\undermod A\ra \overmod A^{\op}$
such that the composition $D\Tr:\undermod A\ra\overmod A$ is an
equivalence with inverse equivalence $\Tr D:\overmod A\ra
\undermod A$. Then, for each indecomposable non-projective
$A$-module $M$, we have $\tau M\cong D\Tr M$, and for each
indecomposable non-injective $A$-module $N$, we have $\tau^{-1}
N\cong \Tr D N$.

For any two $A$-modules $M=\oplus_i M_j$ and $N=\oplus_jN_j$,
where $M_i$ and $N_j$ are indecomposable, we define the {\it
radical} of $\hom_A(M,N)$ by
$$\rad_A(M,N)=\{(f_{ji}):M\ra N\mid f_{ji}:M_i\to N_j
\text{ is not an isomorphism for all }i,j\}.$$ In fact,
$\rad_A(-,-)$ is an ideal of $\mod A$. Inductively, for each
$n>1$, the $n$-th power of the radical is defined to be
$$\rad^n_A(M,N)=\sum_{X}\rad_A^{n-1}(X,N)\circ\rad_A(M,X).$$

From the functorial point of view (see, e.g., \cite[4.8]{Be}), we
may study almost split sequences in the category ${\mathbf
Fun}(A)$ (resp. ${\mathbf Fun}^{\op}(A)$) whose objects are the
covariant (resp. contravariant) additive functors from $\modh A$
(resp. $\mod A^{\op}$) to the category ${\mathbf vec}_k$ of
$k$-vector spaces, and whose morphisms are the natural
transformations of functors. Thus, each $A$-module $M$ defines a
functor $\hom_A(-,M)$ in ${\mathbf Fun}^{\op}(A)$ by
$\hom_A(-,M)(N)=\hom_A(N,M)$ and a subfunctor $\rad_A(-,M)$ by
$\rad_A(-,M)(N)=\rad_A(N,M)$. We denote the induced quotient
functor by
$$\sH_M:=\hom_A(-,M)/\rad_A(-,M).$$
In particular, if $M$ is indecomposable, $\sH_M$ is a simple
functor
--- a functor that has no non-zero proper subfunctor.
Dually, each $A$-module $M$ defines a
functor $\hom_A(M,-)$ in ${\mathbf Fun}(A)$ and the quotient
functor ${\cal K}_M$ by its subfunctor $\rad_A(M,-)$.
We have the
following functorial characterization of minimal right or left
almost split morphisms;
see \cite[1.4]{Ga2} or \cite[4.12.6]{Be}.

\begin{prop} \label{fun}
\begin{itemize}
\item[(a)] A morphism $\vphi:L\ra M$ of $A$-modules is
minimal right almost split if and only if the induced sequence
$$0\lra \hom_A(-,\ker\vphi)\skrel{\kappa_*}{\lra}\hom_A(-,L)
\skrel{\vphi_*}{\lra}\hom_A(-,M)\skrel{\xi_M}{\lra} \sH_M\lra 0$$
is a minimal projective resolution of $\sH_M$ in ${\mathbf
Fun}^{\op}(A)$, where $\kappa:\Ker\vphi\ra L$ is the canonical inclusion
and $\xi_M$ denotes the canonical projection.
\item[(b)] A morphism $\psi:N\ra L$ of $A$-modules is
minimal left almost split if and only if the induced sequence
$$0\lra \hom_A(\Coker\psi,-)\skrel{\pi_*}{\lra}\hom_A(L,-)
\skrel{\psi_*}{\lra}\hom_A(N,-)\skrel{\xi_N}{\lra} {\cal K}_N\lra 0$$
is a minimal projective resolution of ${\cal K}_N$ in ${\mathbf
Fun}(A)$, where $\pi:L\ra \Coker\psi$ is the canonical projection.
\end{itemize}
\end{prop}

{\it We now assume $k=\ofq$} and let $F_A$ be a fixed Frobenius
morphism on $A$. We also assume that each $A$-module $M$ has an
$\BF_q$-structure given by a Frobenius map $F_M$.

\begin{lem} \label{MOR} A morphism
$\vphi: L\ra M$ (resp. $\psi:N\to L$) in $\mod A$ is minimal right
(resp. left) almost split if and only if so is $\vphi^{[1]}:
L^{[1]}\ra M^{[1]}$ (resp. $\psi^{[1]}:N^{[1]}\to L^{[1]}$). In
other words, a sequence $0\ra
N\skrel{\psi}{\ra}L\skrel{\vphi}{\ra}M\ra 0$ is an almost split
sequence if and only if so is $0\ra
N^{[1]}\skrel{\psi^{[1]}}{\ra}L^{[1]}\skrel{\vphi^{[1]}}{\ra}M^{[1]}\ra
0$. Moreover, $\tau M^{[1]}\cong (\tau M)^{[1]}$.
\end{lem}
\begin{pf} Since the Frobenius twisting functor $(\,\,)^{[1]}$ is
invertible with the inverse $(\,\,)^{[-1]}$ given in (\ref{inv}),
everything is clear.
\end{pf}

Let $\vphi:L\ra M$ be a minimal right almost split morphism in
$\mod A$. Then, by the lemma above, for each integer $s\geq 1$,
$\vphi^{[s]}: L^{[s]}\ra M^{[s]}$ is minimal right almost split.
Let $r=p(M)$ be the $F$-period of $M$, i.e., $r$ is minimal with
$M^{[r]}\cong M$. Then Theorem \ref{AR} implies $L^{[r]}\cong L$.
By Proposition \ref{twist}(b), we may assume that $M^{[r]}=M$ and
$L^{[r]}=L$. Thus, both ${\tilde M}=\oplus_{i=0}^{r-1}M^{[i]}$ and
${\tilde L}=\oplus_{i=0}^{r-1}L^{[i]}$ defined in (\ref{Mtil}) are
$F$-stable with respect to Frobenius maps $F_{\tilde M}$ and
$F_{\tilde L}$ defined by $F_M$ and $F_L$, respectively (see
(\ref{eq:tiF})). Since both $\vphi:L\ra M$ and $\vphi^{[r]}: L\ra
M$ are minimal right almost split, there is a $g\in \Aut_A(L)$
such that $\vphi =\vphi^{[r]}\circ g$. Let $F=F_{(L,L)}$ and
$F_{(L,M)}$ be the induced Frobenius maps on $\hom_A(L,L)$
and $\hom_A(L,M)$, respectively; see Lemma \ref{Lang}.
Then $F_{(L,M)}(\psi\circ f)=F_{(L,M)}(\psi)\circ F(f)$ for
all $\psi\in\hom_A(L,M)$ and $f\in\End_A(L)$. Restricting $F$
to the connected algebraic group $\Aut_A(L)$ and applying
Lang-Steinberg's theorem, we may find an element $h\in\Aut_A(L)$
satisfying $g=F^r(h)h^{-1}$, that is, $h=g^{-1}F^r(h)$. This
implies that
$$\vphi \circ h=\vphi\circ(g^{-1}F^r(h))=\vphi^{[r]}\circ F^r(h)
=F_{(L,M)}^r(\vphi\circ h).$$ Note that, by Cor. \ref{twist2},
$F_{(L,M)}^r(\ph)=\ph^{[r]}$. Thus, $\vphi\circ h:L\ra M$ is
minimal right almost split satisfying $(\vphi\circ
h)^{[r]}=\vphi\circ h$. Replacing $\vphi$ by $\vphi\circ h$, we
may assume that $\vphi$ is chosen to satisfy $\vphi^{[r]}=\vphi$.
Now we define
$${\tilde \vphi}:=\diag(\vphi,\vphi^{[1]},\cdots,\vphi^{[r-1]}):
{\tilde L}=\bigoplus_{i=0}^{r-1}L^{[i]}\ra
\bigoplus_{i=0}^{r-1}M^{[i]}=\tilde M.$$ The equality
$\vphi^{[r]}=\vphi$ implies that $\tilde \vphi$ is a morphism in
$\modf A$. Hence, $\tilde\vphi$ induces an $A^F$-module morphism
${\tilde \vphi}^F: {\tilde L}^{F}\ra {\tilde M}^{F}$. Here again
we drop the subscripts of the $F$'s for notational simplicity.

\begin{thm} \label{almor}
Let $\vphi: L\ra M$  be a minimal right almost split morphism.
Then there exists an induced morphism $\tilde \vphi:\tilde L\to
\tilde M$  of $F$-stable modules such that its restriction
${\tilde \vphi}^F: {\tilde L}^{F} \ra {\tilde M}^{F}$  is a
minimal right almost split morphism in $\mod A^F$. In particular,
every almost split sequence $0\ra
N\skrel{\psi}{\ra}L\skrel{\vphi}{\ra}M\ra 0$ in $\modh A$ gives
rise to an almost split sequence
$$0\lra {\tilde N}^{F}\skrel{{\tilde \psi}^F}{\lra}
{\tilde L}^{F}\skrel{{\tilde \vphi}^F}{\lra} {\tilde M}^{F}\lra
0$$ in $\mod A^F$. Moreover, every almost split sequence of
$A^F$-modules can be constructed in this way.
\end{thm}

\begin{pf} Clearly, the last assertion follows from Theorem
\ref{INDF}. From the construction  of $\tilde \vphi$ above, we
have $\Ker{\tilde\vphi}=\oplus_{i=0}^{r-1}\Ker\vphi^{[i]}$ and the
Frobenius map $F_\tL$ on $\tilde L$ induces a Frobenius map on
$\Ker{\tilde\vphi}$. Thus, $\Ker{\tilde\vphi}$ is $F$-stable.

For each $0\leq i\leq r-1$, by Prop. \ref{fun} and Lemma
\ref{MOR}, the morphism $\vphi^{[i]}:L^{[i]}\ra M^{[i]}$ gives the
minimal projective resolution of $\sH_{M^{[i]}}$
$$0\lra \hom_A(-,\Ker\vphi^{[i]})\skrel{\kappa^{[i]}_*}{\lra}\hom_A(-,L^{[i]})
\skrel{\vphi^{[i]}_*}{\lra}\hom_A(-,M^{[i]})\skrel{\xi_{M^{[i]}}}{\lra}
\sH_{M^{[i]}}\lra 0,$$ where $\kappa$ is the inclusion
$\Ker{\tilde\vphi}^{[i]}\to L^{[i]}$. Summing up, we obtain a
minimal projective resolution of $\sH_{\tilde
M}\cong\oplus_{i=0}^{r-1}\sH_{M^{[i]}}$
$$0\lra \hom_A(-,\Ker{\tilde \vphi})\skrel{\tilde \kappa_*}{\lra}
\hom_A(-,\tilde L)\skrel{\tilde \vphi_*}{\lra} \hom_A(-,{\tilde
M})\skrel{\xi_{\tilde M}}{\lra} \sH_{\tilde M}\lra 0,$$ where
$\tilde \kappa=\diag\{\kappa_i\}:\Ker{\tilde
\vphi}=\oplus_{i=0}^{r-1}\Ker\vphi^{[i]} \ra {\tilde L}$. Thus,
for each $F$-stable module $(X,F_X)$, we get the following exact
sequence
$$0\lra \hom_A(X,\Ker{\tilde \vphi})\skrel{\tilde \kappa_*(X)}{\lra}
\hom_A(X,\tilde L)\skrel{\tilde \vphi_*(X)}{\lra} \hom_A(X,{\tilde
M})\skrel{\xi_{\tilde M}(X)}{\lra} \sH_{\tilde M}(X)\lra 0.$$ Now
the Frobenius maps on modules induce Frobenius maps $\tilde F$ on
each space in the sequence above. It is easy to see that all
morphisms in the sequence are compatible with those Frobenius maps
$\tilde F$. This gives the exact sequence
$$0\lra \hom_A(X,\Ker{\tilde \vphi})^{\tilde F}\lra
\hom_A(X,\tilde L)^{\tilde F}\lra\hom_A(X,{\tilde M})^{\tilde F}
\lra \sH_{\tilde M}(X)^{\tilde F}\lra 0.$$ Then we deduce from
Prop. \ref{soc} the following exact sequence
$$0\ra \hom_{A^F}(X^{F},\Ker({\tilde \vphi}^F))\ra
\hom_{A^F}(X^{F},{\tilde L}^F)\ra \hom_{A^F}(X^{F},{\tilde
M}^F)\ra \sH_{{\tilde M}^F}(X^{F})\ra 0.$$ Since every
$A^F$-module is of the form $X^F$ (Thm \ref{CATISO}), we obtain a
minimal projective resolution of the simple functor $\sH_{{\tilde
M}^F}$ in ${\mathbf Fun}(A^F)$
$$0\lra \hom_{A^F}(-,\Ker({\tilde \vphi}^F))\skrel{\tilde\kappa_*^F}\lra
\hom_{A^F}(-,{\tilde L}^F)\skrel{{\tilde \vphi}^F_*}{\lra}
\hom_{A^F}(-,{\tilde M}^F)\skrel{\xi_{\tilde M}^F}{\lra}
\sH_{{\tilde M}^F}\lra 0.$$ Therefore, by Lemma \ref{fun} again,
${\tilde \vphi}^F: {\tilde L}^F \ra {\tilde M}^F$ is a minimal
right almost split morphism in $\mod A^F$.
\end{pf}

Dually, for each minimal left almost split morphism $\psi: N\ra L$
in $\mod A$,  we can construct a minimal left almost split
morphism ${\tilde \psi}^F: {\tilde N}^{F}\ra {\tilde L}^{F}$ in
$\mod A^F$ in a similar way. We leave the detail to the reader.

\section{The Auslander-Reiten quivers}

We are now going to prove that the Frobenius morphism $F$ on $A$
induces an automorphism $\fs$ of the Auslander-Reiten quiver $\cq$
of $A$ and that the induced modulated quiver $\fM_{\cq,\fs}$ is
essentially the Auslander-Reiten quiver of $A^F$.

We begin with the general definition. Let $A$ be a finite
dimensional algebra over an {\it arbitrary} field $k$. For an
$A$-module $M$, let $D_M$ denote the $k$-algebra
$$D_M:=\End_A(M)/\rad(\End_A(M)).$$
This is a division algebra if $M$ is indecomposable. By
definition, the {\it Auslander-Reiten quiver} (or {\it AR-quiver}
for short) of $A$ is a (simple) $k$-modulated quiver ${\cal Q}_A$
consisting of a valued graph $\Ga=\Ga_A$ 
and a $k$-modulation $\BM=\BM_A$ defined on
$\Ga$. Here, the vertices of $\Ga$ are isoclasses $[M]$ of
indecomposable $A$-modules and the arrows $[M]\ra [N]$ for
indecomposable $M$ and $N$ are defined by the condition $\Irr_A(M,
N)\not=0$, where
$$\Irr_A(M,N):=\rad_A(M,N)/\rad_A^2(M,N)$$
is the space of irreducible homomorphisms from $M$ to $N$. Each
arrow $[M]\ra [N]$ has the valuation $(d_{MN},d'_{MN})$ with
$d_{MN}$ and $d'_{MN}$ being the dimensions of $\Irr_A(M,N)$
considered as left $D_N$-space and right $D_M$-space,
respectively. The $k$-modulation $\BM$ is given by division
algebras $D_M$ for vertices $[M]$ and (non-zero)
$D_N$-$D_M$-bimodules $\Irr_A(M,N)$ for arrows $[M]\ra [N]$.

\begin{rems} \label{AR-quiver} (a) The AR-quiver of $A$ defined
in  \cite[VII.1]{ARS} is simply the valued quiver $\Ga_A$ together
with the translation $\tau$ sending non-projective vertices to
non-injective vertices. The modulation $\BM_A$ for ${\cal Q}_A$ is
not explicitly mentioned there. Here we adopt the definition for
AR-quivers given by Benson in \cite[p.150]{Be}. We will see that
this definition is more natural to fit our situation. Moreover, we
shall also see that the translation $\tau$ for $A^F$ is naturally
induced from the translation for $A$.

(b) The valuation $(d,d')$ of an arrow $[M]\ra[N]$ in $\cq_A$
admits the following description (see \cite[VII.1]{ARS}). If $L\ra
N$ is minimal right almost split, then $L\cong d'M\oplus L_1$,
where $L_1$ has no summand isomorphic to $M$. If $M\ra K$ is
minimal left almost split, then $K\cong dN\oplus K_1$, where $K_1$
admits no summand isomorphic to $N$.

(c) If $k$ is algebraically closed, then $D_M\cong k$ and
$d_{MN}=d'_{MN}=\dim_k\Irr(M,N)$ for all indecomposable
$A$-modules $M$ and $N$. So the modulation in the AR-quiver
$\cq_A$ consists of $k$-spaces which can be represented by drawing
$d_{MN}$ arrows from $[M]$ to $[N]$. In this way, we turn the
modulated quiver $\cq_A$ to an ordinary quiver.
\end{rems}

We return to our usual setup: let $k=\ofq$ and let $F_A$ be a
fixed Frobenius morphism on $A$. For each indecomposable
$A$-module $M$, let $r=p(M)$ be the period of $M$. In view of
Prop. \ref{twist}(b), we may assume that $M^{[r]}=M$. Then
$\tM=M\oplus M^{[1]}\oplus \cdots\oplus M^{[r-1]}$ is $F$-stable
with $F_\tM:\tM\ra\tM$ defined by
$$F_\tM(x_0,x_1,\cdots,x_{r-1})
=(F_M(x_{r-1}),F_M(x_0),\cdots,F_M(x_{r-2})).$$ By Theorem
\ref{INDF}, $\tM^{F}$ is an indecomposable $A^F$-module. By Prop.
\ref{soc}, the induced Frobenius map $F=F_{(\tM,\tM)}$ on
$\End_A(\tM)$ gives a canonical $\bbf_q$-algebra isomorphism
$$\bigl(\End_A(\tM)/\rad\End_A(\tM)\bigr)^F\cong
\End_{A^F}(\tM^{F})/\rad\End_{A^F}(\tM^{F}).$$ In particular, if
$\tM=X_k:=X\otimes k$ for an indecomposable $A^F$-module $X$, then
$$\bigl(\End_A(\tM)/\rad\End_A(\tM)\bigr)^F=(D_{X_k})^F\cong D_X.$$

Let $X$ and $Y$ be indecomposable $A^F$-modules. Up to
isomorphism, we may assume that $X=\tM^{F}$ and $Y=\tN^{F}$ for
some indecomposable $A$-modules $M$ and $N$ (see the proof of Thm
\ref{INDF}). In fact, we may choose $\tM=X_k$ and $\tN=Y_k$ so
that $M$ and $N$ are indecomposable direct summands of the
$A$-modules $X_k$ and $Y_k$, respectively. We then have
$$D_X\cong(D_\tM)^{F}\cong \bbf_{q^{r_1}}\;\;\text{and}\;\;
D_Y\cong(D_\tN)^{F}\cong \bbf_{q^{r_2}},$$ where $r_1=p(M)$ and
$r_2=p(N)$. Moreover, the Frobenius map $F_{(\tM,\tN)}$ on
$\hom_A(\tM,\tN)$ induces Frobenius maps on $\rad_A^n(\tM,\tN)$
for each $n\geq 1$, and thus, a Frobenius map $F$ on
$\Irr_A(\tM,\tN)$. By viewing $\bigl(\Irr_A(\tM,\tN)\bigr)^F$ as a
$D_Y$-$D_X$-bimodules via the isomorphisms $(D_\tM)^{F}\cong D_X$
and $(D_\tN)^{F}\cong D_Y$, we have the following lemma.

\begin{lem} \label{modulation} Let $X=\tM^F$ and $Y=\tN^F$ be indecomposable
$A^F$-modules, where $M$ and $N$ are indecomposable $A$-modules.
 Then the $D_Y$-$D_X$-bimodules $\Irr_{A^F}(X,Y)$ and $\bigl(\Irr_A(\tM,\tN)\big)^F$ are
isomorphic.
\end{lem}

\begin{pf} It is clear to see that restriction of the linear isomorphism
given in Prop. \ref{soc}(b) gives an $\bbf_q$-linear isomorphism
$$\Psi: \rad_{A^F}(X,Y)\overset\sim\lra \bigl(\rad_A(\tM,\tN)\bigr)^{F};
f\lmto f\otimes 1,$$ and hence, an $\bbf_q$-linear injection
$$\Psi:\rad_{A^F}^2(X,Y)\ra\bigl(\rad_A^2(\tM,\tN)\bigr)^{F}.$$
Thus, $\Psi$ induces a surjective map
$${\bar \Psi}: \Irr_{A^F}(X,Y)\ra
\bigl(\rad_A(\tM,\tN)\bigr)^{F}/\bigl(\rad_A^2(\tM,\tN)\bigr)^{F}
\cong \bigl(\Irr_A(\tM,\tN)\big)^F,$$ which is clearly a
$D_Y$-$D_X$-bimodule homomorphism. We now prove that $\bar\Psi$ is
a linear isomorphism by a comparison of dimensions.

For each $0\leq s\leq p(M)-1$ and each $0\leq t\leq p(N)-1$, we
set $n_{st}=\dim_k\Irr_A(M^{[s]},N^{[t]})$. Since
$$\Irr_A(\tM,\tN)\cong \bigoplus_{{0\leq s\leq p(M)-1}\atop{0\leq t\leq p(N)-1}}
\Irr_A(M^{[s]},N^{[t]}),$$ we have
$n:=\dim_k\Irr_A(\tM,\tN)=\sum_{s,t}n_{st}$, and so
$\dim_{\bbf_q}\bigl(\Irr_A(\tM,\tN)\big)^F=n$. Now, take a minimal
 right almost split map $L\to N$. Then, so is $L^{[t]}\to N^{[t]}$
by Lemma \ref{MOR}. For each fixed $s$, it holds that
$L^{[t]}\cong n_{st}M^{[s]}\oplus L_{s,t}$ with $M^{[s]}\nmid
L_{s,t}$ for all $0\le t\le p(N)-1$. Thus, $\tL\cong
n_sM^{[s]}\oplus L_s$, for all $s$ with $0\le s\le p(M)-1$, where
$n_s=\sum_{t=0}^{p(N)-1}n_{st}$ and
$L_s=\oplus_{t=0}^{p(N)-1}L_{s,t}$. Since $\tL=\tL^{[1]}\cong
n_sM^{[s+1]}\oplus L_s^{[1]}$ and the $M^{[s]}$ are pairwise
non-isomorphic, it follows that $n_s=n_{s+1}$ for all $0\leq s<
p(M)-1$. Thus, $n_s=\frac n{p(M)}$ and
$\tL=\frac{n}{p(M)}\tM\oplus\tL'$ for some $F$-stable module
$\tL'$ with $\tM\nmid \tL'$.

On the other hand, by Theorem \ref{almor} (see also Remark
\ref{AR-quiver}(b)), $\tL^F\to\tN^F=Y$ is minimal right almost
split and
$$\tL^F=\frac{n}{p(M)}\tM^F\oplus\tL^{\prime F}=\frac{n}{p(M)}X\oplus\tL^{\prime F}$$
with $X\nmid \tL^{\prime F}$, it follows that, if $(d,d')$ be the
valuation of the arrow $[X]\ra[Y]$ in the AR-quiver $\cq_{A^F}$ of
$A^F$, then $d'=\frac{n}{p(M)}$. Hence,
$$\dim_{\bbf_q} \Irr_{A^F}(X,Y)=d'\dim_{\bbf_q}D_X=n=
\dim_{\bbf_q}\bigl(\Irr_A(\tM,\tN)\big)^F.$$ Consequently, ${\bar
\Psi}$ is a $D_Y$-$D_X$-bimodule isomorphism.
\end{pf}

Since the algebra $A$ is defined over the algebraically close
field $k=\ofq$, we may regard the AR-quiver $\cq=\cq_A$ of $A$ as
an ordinary quiver (see Remark \ref{AR-quiver}(c)). We first
observe that $\cq$ admits an admissible automorphism $\fs$. For
each vertex $[M]\in\cq$, $\fs([M])$ is defined to be $[M^{[1]}]$.
If $M$ and $N$ are indecomposable $A$-modules, then there are
$n_{st}$ arrows $\gz_{s,t}^{(m)}$ from $[M^{[s]}]$ to $[N^{[t]}]$
in $\cq$, where $0\leq s\leq p(M)-1$, $0\leq t\leq p(N)-1$,
$n_{st}=\dim_k\Irr_A(M^{[s]},N^{[t]})$ and $1\leq m\leq n_{st}$.
Note that $n_{st}=n_{s+1,t+1}$ for all $s,t$, where subscripts are
considered as integers modulo $p(M)$ and $p(N)$, respectively. We
now define
$$\fs(\gz_{s,t}^{(m)})=\gz_{s+1,t+1}^{(m)}\;\;\text{for all
$0\leq s\leq p(M)-1$ and $0\leq t\leq p(N)-1$}.$$ Clearly, $\fs$
is an admissible quiver automorphism and $(\cq,\fs)$ is an
ad-quiver.

Associated to $(\cq,\fs)$, we may define a modulated quiver
$\fM_{\cq,\fs}$ as in (\ref{mqs}): let $\ca=k\cq$ denote the path
algebra of $\cq$ and $F=F_{\cq,\fs}$ be the Frobenius morphism of
$\ca$ induced by the automorphism $\fs$. For each vertex $\bi(M)$
(i.e., the $\fs$-orbit of $[M]$) and each arrow $\brho$ (i.e., an
$\fs$-orbit of arrows in $\cq$) in $\ggz(\cq,\fs)$, we define
subspaces
$$\ca_{\bi(M)}=\bigoplus_{s=0}^{p(M)-1}ke_{[M^{[s]}]}\;\;\text{and}\;\;
\ca_\brho=\bigoplus_{\rho\in\brho}k\rho,$$ of $\ca$, which are
obviously $F$-stable. By definition, the $\bbf_q$-modulation
$\BM({\cq,\fs})$ is given by $(\ca_{\bi(M)})^{F}$ and
$(\ca_\brho)^{F}$ for all vertices $\bi(M)$ and arrows $\brho$ in
$\ggz(\cq,\fs)$.

Recall the definition \ref{mqi} of isomorphisms for modulated
quivers. We now can state the following result.

\begin{thm} \label{ARqiso} The modulated quiver $\fM_{\cq,\fs}$
associated to the AR-quiver $(\cq,\fs)$ of $A$ defined above is
isomorphic to the AR-quiver $\cq_{A^F}$ of $A^F$. Moreover, the
Auslander-Reiten translation of $A$ naturally induces that of the
fixed-point algebra $A^F$.
\end{thm}

\begin{pf} If $X=\tM^{F}$ is an indecomposable $A^F$-module, then clearly,
the correspondence $[X]\lra \bi(M)$ gives a bijection between
vertices of $\cq_{A^F}$ and those of $\fM_{\cq,\fs}$ (see
Thm\ref{INDF}). Moreover, there is an isomorphism
\begin{equation}\label{DtM}
\ca_{\bi(M)}=\bigoplus_{s=0}^{p(M)-1}ke_{[M^{[s]}]}\lra
\bigoplus_{s=0}^{p(M)-1}D_{M^{[s]}}\cong D_\tM,\; \sum_{s}x_s
e_{[M^{[s]}]}\lmto (x_s1_{M^{[s]}})_s
\end{equation} which induces
an $\bbf_q$-algebra isomorphism $(\ca_{\bi(M)})^{F}\cong
(D_{\tM})^F\cong D_X$. As before, we shall identify
$(\ca_{\bi(M)})^{F}$ with $D_X$. Thus, $(\ca_\brho)^{F}$ is a
$D_Y$-$D_X$-bimodule.

It remains to prove that for $X=\tM^{F}$ and $Y=\tN^{F}$, where
$M$ and $N$ are indecomposable $A$-modules, there is a
$D_Y$-$D_X$-bimodule isomorphism
\begin{equation}\label{bimod}
\bigoplus_{\brho:\bi(M)\ra \bi(N)}(\ca_{\brho})^{F}\cong
\Irr_{A^F}(X,Y).
\end{equation}
 By Lemma \ref{modulation}, it suffices to show that we
have a $(D_{\tN})^F$-$(D_{\tM})^F$-bimodule isomorphism
$$\bigoplus_{\brho:\bi(M)\ra \bi(N)}(\ca_\brho)^{F}\cong
\bigl(\Irr_A(\tM,\tN)\bigr)^F,$$ or a $D_{\tN}$-$D_{\tM}$-bimodule
isomorphism
\begin{equation}\label{*}
\phi:\bigoplus_{\brho:\bi(M)\ra \bi(N)}\ca_\brho\overset\sim\lra
\Irr_A(\tM,\tN) \end{equation} which is compatible with the
Frobenius morphism $F_\ca$ on $\ca$ and $F$ on $\Irr_A(\tM,\tN)$.

 Let $d$ and $l$ denote the
greatest common divisor and the least common multiple of
$p_1=p(M)$ and $p_2=p(N)$, respectively. As before, let
$n_{st}=\dim_k\Irr_A(M^{[s]},N^{[t]})$ and denote by
$\gz_{s,t}^{(m)}$ all arrows from $[M^{[s]}]$ to $[N^{[t]}]$ in
$\Ga(\cq,\fs)$. Then the set of arrows
$$\{\gz_{s,t}^{(m)}|0\leq s\leq p_1-1, 0\leq t\leq p_2-1,
1\leq m\leq n_{st}\}$$ is $\fs$-stable and the arrows
$\gz_{0,t}^{(m)}$ with $0\leq t\leq d-1$ and $1\leq m\leq n_{0t}$
form a complete set of representatives of $\fs$-orbits in this
set, each of which is of length $l$. On the other hand, for each
$0\leq t\leq d-1$, the space $\Irr_A(M,N^{[t]})$ is $F^l$-stable
with $l$ minimal. Thus, we can choose a $k$-basis
$\xi^{(m)}_{0,t}$, $1\leq m\leq n_{0t}$, of $\Irr_A(M,N^{[t]})$
which are $F^l$-fixed. Then the set
$$\{F^a(\xi_{0,t}^{(m)})\mid 0\leq t\leq d-1,1\leq m\leq n_{0t},
0\leq a\leq l-1\}$$ is a $k$-basis for
$\bigoplus_{s,t}\Irr_A(M^{[s]},M^{[t]})\cong \Irr_A(\tM,\tN)$.
Thus, the correspondence $\gz_{0,t}^{(m)}\lmto \xi_{0,t}^{(m)}$
induces an isomorphism of $k$-spaces
$$\phi:\bigoplus_{\brho:\bi(M)\ra \bi(N)}\ca_\brho
=\bigoplus_{s,t,m}k\gz_{s,t}^{(m)}\overset\sim\lra
\Irr_A(\tM,\tN)$$ sending $\ga_{s,t+s}^{(m)}$ to
$F^s(\xi_{0,t}^{(m)})$. Therefore, $\phi\circ F_\ca=F\circ \phi$.
Using the isomorphism given in (\ref{DtM}), we see easily that
$\phi$ is a $D_{\tN}$-$D_{\tM}$-bimodule isomorphism. So (\ref{*})
is proved.

To prove the last assertion, we first note that the Frobenius
morphism $F_A$ on $A$ is also a Frobenius morphism on $A^{\op}$.
Since the projective cover of an $A$-module is unique up to
isomorphism, each $F$-stable $A$-module $M$ admits a minimal
projective presentation $P_1\skrel{f}{\lra}P_0\skrel{\pi}{\lra
M}\lra 0$ such that $P_0$ and $P_1$ are $F$-stable projective and
that $f$ and $\pi$ are morphisms in $\modf A$. Thus the transpose
$\Tr: \undermod A\ra \undermod A^{\op}$ induces a functor
$\ct:\undermodf A\ra \undermodf A^{\op}$, where, for example,
$\undermodf A$ denotes the factor category of $\modf A$ modulo
morphisms factoring through an $F$-stable projective $A$-module.
It is obvious that the dual $D=\hom_k(-,k):\mod A\ra\mod A^{\op}$
induces a dual $\cd:\undermodf A\ra\overmodf A^{\op}$. On the
other hand, for the fixed algebras $A^F$ and
${(A^{\op})}^F={(A^F)}^{\op}$, we also have the transpose
$\Tr^F:\undermod A^F\ra \undermod {(A^F)}^{\op}$ and the dual
$D^F=\hom_{\bbf_q}(-,\bbf_q):\undermod A^F\ra\overmod
{(A^F)}^{\op}$. Furthermore, the functor $\Phi=\Phi_A:\modf A\ra
\mod A^F$ defined in the proof of Theorem \ref{CATISO} induces
functors $\Phi_A:\undermodf A\ra \undermod A^F$ and
$\Phi_A:\overmodf A\ra \overmod A^F$. By Proposition \ref{soc}, we
finally have the following commutative square
$$\unitlength=1cm
\begin{picture}(6,3.2)
\put(0.2,2.5){$\undermodf A$} \put(4,2.5){$\overmodf A^{\op}$}
\put(0.2,0.5){$\undermod A^F$} \put(4,0.5){$\overmod
{(A^{\op})}^F$} \put(1.9,2.6){\vector(1,0){1.8}}
\put(1.9,0.6){\vector(1,0){1.8}}
\put(0.75,2.3){\vector(0,-1){1.4}}
\put(4.75,2.3){\vector(0,-1){1.4}} \put(2.5,2.7){$\cd\ct$}
\put(2.35,0.7){$D^F\Tr^F$} \put(0.2,1.6){$\Phi_A$}
\put(5,1.6){$\Phi_{A^{\op}}$}
\end{picture}$$
\end{pf}

\begin{rems} (1) By the definition of $\fs$, each $(\ca_\brho)^F$ is
a simple $D_Y$-$D_X$-bimodule. Thus, the isomorphism in
(\ref{bimod}) gives a decomposition of the $D_Y$-$D_X$-bimodule
$\Irr_{A^F}(X,Y)$ into the sum of simples.

(2) By Thms \ref{almor} and \ref{ARqiso}, we see that the
Auslander-Reiten theory for algebras defined over a finite field
$\BF_q$ is completely determined by the theory for algebras
defined over the algebraic closure ${\ol \BF}_q$.
\end{rems}

\section{Counting the number of $F$-stable representations}

In this section, we shall use the theory established in \S 6 to
count the number of representations (resp. indecomposable
representations) of a finite dimensional hereditary (basic)
$\BF_q$-algebra in terms of representations of an ad-quiver. In
particular, we shall prove that, for a given dimension vector,
these numbers are polynomials in $q$.

Let $(Q,\si)$ be an ad-quiver and $\fM_{Q,\si}$ the associated
modulated quiver with underlying valued quiver
$\ggz=\ggz(Q,\sz)=(I,\ggz_1)$. By definition, a representation
$V=(V, \phi)$ of $Q$ over $ k$ consists of a $Q_0$-graded
$k$-vector space $V=\oplus_{i\in Q_0}V_i$ and a family of $
k$-linear maps $\phi=(\phi_\rho)_\rho$ with $\phi_\rho:V_i\ra V_j$
for all arrows $\rho:i\ra j\in Q_1$. A morphism from $(V,\phi)$ to
$(V',\phi')$ is given by a $Q_0$-graded morphism $f=(f_i)_i:V\ra
V'$ such that $f_{t\rho}\circ \phi_\rho=\phi_\rho'\circ f_{h\rho}$
for each arrow $\rho$. We denote by $\rep Q=\rep_kQ$ the category
of all finite dimensional representations of $Q$ over $ k$. It is
well known that $\mod A$, where $A=kQ$, and $\rep Q$ are
isomorphic categories. For each $Q_0$-graded vector space
$V=\oplus_{i\in Q_0}V_i$, we set
$$\udim V=\sum_{i\in I}(\dim_ k V_i)i\in\bbn Q_0,$$
called the {\it dimension vector} of $V$.

Following Lusztig \cite{L98}, by $\cc_I$ we denote the
sub-category of the category $\sV_{k,\BF_q}$ (see the definition
before \ref{Lang}) whose objects are $Q_0$-graded $k$-vector
spaces $V=\oplus_{i\in Q_0}V_i$ with an $\bbf_q$-rational
structure determined by a Frobenius map $F=F_V:V\ra V$ such that
$F(V_i)=V_{\sz(i)}$ for each $i\in Q_0$. The morphisms in $\cc_I$
are $ k$-linear maps compatible with both the gradings and the
Frobenius maps. Note that, for each $\az=\sum_{i\in
Q_0}a_ii\in\bbn Q_0$ satisfying $a_i=a_{\sz (i)}$, $\forall i\in
Q_0$, there is a unique, up to isomorphism, object in $\cc_I$ with
dimension vector $\az$ (see remarks before \ref{Lang}).

A representation $(V,\phi)$ of $Q$ with $V\in\cc_I$ is called
$F$-{\it stable} if, for each $\az\in Q_1$, the equality
$F\circ\phi_\rho=\phi_{\sz(\rho)}\circ F$ holds, that is, the
following diagram commutes
$$\unitlength=1cm
\begin{picture}(5,3)
\put(0.7,2.5){$V_{t\rho}$} \put(4,2.5){$V_{h\rho}$}
\put(0.7,0.5){$V_{\sz(t\rho)}$} \put(4,0.5){$V_{\sz(h\rho)}$}
\put(1.55,2.6){\vector(1,0){2.2}}
\put(1.55,0.75){\vector(1,0){2.2}}
\put(0.9,2.4){\vector(0,-1){1.5}}
\put(4.2,2.4){\vector(0,-1){1.5}} \put(2.5,2.8){$\phi_\rho$}
\put(2.5,0.95){$\phi_{\sz(\rho)}$} \put(0.35,1.7){$F$}
\put(4.3,1.7){$F$}
\end{picture}$$

Let $\rep^F(Q,\sz)$ be the category consisting of $F$-stable
representations of $Q$ together with morphisms in $\rep Q$ which
are compatible with Frobenius maps. Clearly, $\rep^F(Q,\sz)$ is a
subcategory of $\rep Q$, though not necessarily full. In fact,
$\rep^F(Q,\sz)$ is an abelian $\bbf_q$-category. It is easy to see
that the category isomorphism between $\mod A$ and $\rep Q$
induces a category isomorphism between $\modf A$, where
$F=F_{Q,\si}$, and $\rep^F(Q,\sz)$. Hence, by Theorem
\ref{CATISO}, categories $\rep^F(Q,\sz)$ and $\modh A^F$ are
isomorphic.

The quiver automorphism $\si$ extends linearly to a group
automorphism $\si$ on $\BZ Q_0$ defined by $\sz(\sum_{i\in
Q_0}a_ii)=\sum_{i\in Q_0}a_i\sz(i)$. Clearly, if $V$ is
$F$-stable, then $\si(\udim V)=\udim V$. Let $(\bbz Q_0)^\sz$
denote the $\si$-fixed-point subset of $\BZ Q_0$. This set can be
identified with the group $\BZ I$ via the canonical bijection
\begin{equation}\label{hatsi}
\hsz:(\bbz Q_0)^\sz\ra \bbz I;\;\; \sum_{i\in Q_0}b_i i\mapsto
\sum_{\bi\in I}a_\bi\bi,\end{equation} where $a_\bi:=b_i=b_j$ for
all $i,j\in\bi$ (see \cite[\S2]{Hub1}). In particular, the {\it
dimension vector} $\udim X=\sum_{\bi\in I} d_\bi\bi\in\BN I$ of an
$A^F$-module $X$ can be defined by
$$\udim X=\hsz(\udim(X\otimes k)).$$
Note that if, for each $\bi\in I$, $S_\bi$ denotes the simple
$A^F$-module corresponding to the idempotent
$e_\bi=\sum_{i\in\bi}e_i$ of $A^F$, then $d_\bi=\dim_{\fqi}e_\bi
X$ is  the number of composition factors isomorphic to $S_\bi$ in
a composition series of $X$.

Given a matrix $x=(x_{ij})\in k^{m\times n}$ and an integer $r\geq 0$, we define
$$x^{[r]}=(x_{ij}^{q^r})\in k^{m\times n}.$$

For each $\bz=\sum_{i\in Q_0}b_ii\in (\bbz Q_0)^\sz$, let $V_i=k^{b_i}$
for each $i\in Q_0$. We define a Frobenius map $F$ on the $Q_0$-graded vector
space $V=\oplus_{i\in Q_0}V_i$ such that, for $v\in V_i$,
$F(v)=v^{[1]}\in V_{\sz(i)}$ for all $i\in Q_0$. We further define
$$R(\bz)=R(Q,\bz)=\prod_{\rho\in Q_1}\hom_k(k^{b_{t\rho}},k^{b_{h\rho}})\cong
\prod_{\rho\in Q_1} k^{b_{h\rho}\times b_{t\rho}}.$$ Then the
Frobenius map $F$ on $V$ induces a Frobenius map on the variety
$R(\bz)$ such that, for $x=(x_\rho)\in R(\bz)$, $F(x)=(y_\rho)$ is
defined by
$$y_\rho(F(v))=F(x_{\si^{-1}(\rho)}(v))
\text{\quad for all \quad}\rho\in Q_1, v\in V_{\sz^{-1}(t\rho)}.$$
By viewing each $x_\rho$ as a matrix over $k$, we have
$y_{\sz(\rho)}=x_\rho^{[1]}$.
Obviously, a point $x=(x_\rho)_\rho$ of $R(\bz)$ determines a representation
$V(x)=(V_i,x_\rho)$ of $Q$. The algebraic subgroup
$$G(\bz)=\prod_{i\in Q_0} GL_{b_i}(k)$$
of $GL(V)$ acts on $R(\bz)$  by conjugation
$$(g_i)_i\cdot(x_\rho)_\rho=(g_{h\rho}x_\rho g_{t\rho}^{-1})_\rho,$$
and the $G(\bz)$-orbits $\co_x$ in $R(\bz)$ correspond
bijectively to the isoclasses $[V(x)]$ of representations of $Q$
with dimension vector $\bz$.

The Frobenius map $F$ on $V$ also induces a Frobenius map on the
group $GL(V)$ given by
$$F(gv)=F(g)F(v)\text{\quad for all \quad}g\in GL(V), v\in V.$$
It is clear that the subgroup $G(\bz)$ is $F$-stable such that
$F(g)_{\sz(i)}=g_i^{[1]}$ for $g=(g_i)\in G(\bz)$. The action of
$G(\bz)$ on $R(\bz)$ restricts to an action of $G(\bz)^F$ on
$R(\bz)^F$. Then, the $G(\bz)^F$-orbits in $R(\bz)^F$ correspond
bijectively to the isoclasses of $F$-stable representations in
$\rep^F(Q,\sz)$ with dimension vector $\hsz(\bz)$, or
equivalently, to the isoclasses of $A^F$-modules with dimension
vector $\hsz(\bz)$.

For each $\az=\sum_{\bi\in I}a_\bi\bi\in \bbn I$, let
$M_{Q,\sz}(\az,q)$ (resp. $I_{Q,\sz}(\az,q)$) be the number of
isoclasses of $A^F$-modules (resp. indecomposable $A^F$-modules)
of dimension vector $\az$.
Further, as indicated above, $M_{Q,\sz}(\az,q)$ is the number of
$G(\bz)^F$-orbits in $R(\bz)^F$, where $\bz=\hsz^{-1}(\az)\in\bbn
Q_0$.

Let $\bz=\sum_{i\in Q_0}b_ii$. Note that $b_i=a_\bi$ for all $i\in\bi$.
For each $\bi\in I$ and each arrow $\brho:\bi\ra\bj$ in $\ggz$, we define
$$G_\bi=\prod_{i\in\bi}GL_{b_i}(k)
\;\;\text{and}\;\; R_\brho=\prod_{\rho\in\brho}\hom_k(k^{b_{t\rho}},k^{b_{h\rho}})
\cong\prod_{\rho\in \brho} k^{b_{h\rho}\times b_{t\rho}}.$$
By fixing an $i_0\in\bi$ and a $\rho_0\in\brho$, we can identify
$$G_\bi^F=\{(g_i)_{i\in\bi}|g_{\sz(i)}=g_i^{[1]} \;\;\text{for all $i\in\bi$}\}$$
with $GL_{a_\bi}(\fqi)$ by $g_\bi:=(g_i)_{i\in\bi}\lmto g_{i_0}$
and identify
$$R_\brho^F=\{(x_\rho)_{\rho\in\brho}|x_{\sz(\rho)}=x_\rho^{[1]} \;\;
\text{for all $\rho\in\brho$}\}$$ with
$\bbf_{q^{\ez_\brho}}^{a_{\bj}\times a_{\bi}}$ by
$x_\brho:=(x_\rho)_{\rho\in\brho}\lmto x_{\rho_0}$. Hence, we may
identify their products
$$G(\bz)^F=\prod_{\bi\in I}G_\bi^F
\;\;\text{and}\;\; R(\bz)^F=\prod_{\brho\in \ggz_1}R_\brho^F,$$
with
$$G:=\prod_{\bi\in I}GL_{a_\bi}(\fqi)
\;\;\text{and}\;\;X:= \prod_{\brho\in
\ggz_1}\bbf_{q^{\ez_\brho}}^{a_{\bj}\times a_{\bi}}$$
 respectively. Under this identification, the action of
$G(\bz)^F$ on $R(\bz)^F$ becomes the action of $G$ on $X$: for
$g=(g_\bi)\in G$ and $x=(x_\brho)\in X$,
$$(g\cdot x)_\brho=g_\bj^{[s_\brho]}x_{\brho}
(g_\bi^{[r_\brho]})^{-1}\;\text{for each arrow
$\brho:\bi\ra\bj$},$$ where $0\leq r_\brho\leq \ez_\bi-1$ and
$0\leq s_\brho\leq \ez_\bj-1$ are determined by
$\rho_0:\sz^{r_\brho}(i_0)\ra\sz^{s_\brho}(j_0)$ (see Remark
\ref{Hua-Hub}(a)). Then $M_{Q,\sz}(\az,q)$ is the number of
$G$-orbits in $X$.

For each $g=(g_\bi)\in G$, we set $X^g=\{x\in X\mid g\cdot x=x\}$
and $G_g=\{h\in  G\mid hg=gh\}$. By Burnside's formula, we have
$$M_{Q,\sz}(\az,q)=\frac1{|G|}\sum_{g\in G}X^g=\sum_{g\in\text{ccl}(G)}\frac{|X^g|}{|G_g|},$$
where ccl$(G)$ is a set of representatives of conjugacy classes of
$G$. Further, we have
$$|X^g|=\prod_{\brho:\bi\ra\bj}|X^g_\brho|\;\;\text{and}\;\;
|G_g|=\prod_{\bi\in I} |GL_{a_\bi}(\fqi)_{g_\bi}|$$ where
$X^g_\brho=\{x_\brho\in \bbf_{q^{\ez_\brho}}^{a_{\bj}\times
a_{\bi}}\mid g_\bj^{[s_\brho]}x_\brho=x_\brho g_\bi^{[r_\brho]}\}$
and $GL_{a_\bi}(\fqi)_{g_\bi}=\{h_\bi\in GL_{a_\bi}(\fqi)\mid
h_\bi g_\bi=g_\bi h_\bi\}$.

In order to compute $M_{Q,\sz}(\az,q)$, we need to deal with the
conjugacy classes in $GL_m(\bbf_{q^r})$ for $m,r\geq 1$ (see, for
example, \cite[Chap.IV]{M}).

We denote by $\Phi(q^r)$ the set all irreducible polynomials
in $T$ over $\bbf_{q^r}$ with leading coefficients $1$, excluding the
polynomial $T$, and by $\cp$ the set of all partitions, i.e. finite sequences
$\lz=(\lz_1, \lz_2,\cdots)$ of non-negative integers with
$\lz_1\geq \lz_2\geq \cdots$. For $\lz,\mu\in\cp$, we define $|\lz|=\sum_i\lz_i$,
$b_\lz(q)=\prod_{i\geq 1}(1-q)(1-q^2)\cdots (1-q^{\lz_i})$,
and $\lr{\lz,\mu}=\sum_{i,j}\text{min}\{\lz_i,\mu_j\}$.

It is known that the conjugacy classes in $GL_m(\bbf_{q^r})$ are in one-to-one
correspondence with the isoclasses of $\bbf_{q^r}[T]$-modules of dimension
$m$. The latter is parametrized by the functions
$\pi:\Phi(q^r)\ra\cp$ such that $\sum_{\vphi\in\Phi(q^r)} d(\vphi)|\pi(\vphi)|=m$,
where $d(\vphi)$ denotes the degree of $\vphi$. Each such a function $\pi$
corresponds to the module
$$\bigoplus_{\vphi\in\Phi(q^r)}\bigoplus_{i\geq 1}\bbf_{q^r}[T]/(\vphi^{\pi_i(\vphi)}),$$
where $\pi(\vphi)=(\pi_1(\vphi),\pi_2(\vphi),\cdots)$.

Let $g\in GL_m(\bbf_{q^r})$ be in the conjugacy class
corresponding to the partition function $\pi:\Phi(q^r)\ra\cp$.
Then we have (see \cite[p.272]{M})
$$|GL_m(\bbf_{q^r})_g|=|\{h\in  GL_m(\bbf_{q^r})|hg=gh\}|=\prod_{\vphi\in\Phi(q^r)}
q^{rd(\vphi)\lr{\pi(\vphi),\,\pi(\vphi)}}b_{\pi(\vphi)}(q^{-rd(\vphi)}).$$

Further, let $s,t,n_1,n_2,n\geq 1$ be such that $n_1|n$ and $n_2|n$. Take
$g_1\in GL_s(\bbf_{q^{n_1}})$ and $g_2\in GL_t(\bbf_{q^{n_2}})$ such that
their conjugacy classes correspond respectively to
partition functions $\pi^1: \Phi(q^{n_1})\ra\cp$ and
$\pi^2: \Phi(q^{n_2})\ra\cp$. By \cite[p.253]{HL}, we have
$$|\{x\in\bbf_{q^n}^{t\times s}|g_2x=xg_1\}|
=q^{n\sum_{\vphi\in\Phi(q^{n_1}), \psi\in\Phi(q^{n_2})}
d(\vphi,\psi)\lr{\pi^1(\vphi),\,\pi^2(\psi)}},$$
where $d(\vphi,\psi)$ denotes the degree of the greatest common divisor
of $\vphi$ and $\psi$ over a common extension field of $\bbf_{q^{n_1}}$ and
$\bbf_{q^{n_2}}$ (Note that $d(\vphi,\psi)$ is independent of the
extension field).

It follows that the conjugacy classes in $G=\prod_{\bi\in I}GL_{a_\bi}(\fqi)$
are in one-to-one correspondence to multi-partition functions
$\pi=(\pi^\bi)_{\bi\in I}$ of $\pi^\bi:\Phi(q^{\ez_\bi})\ra\cp$
with $\sum_{\vphi\in\Phi(q^{\ez_\bi})}d(\vphi)\pi^\bi(\vphi)=a_\bi$. The set
of all such multi-partition functions is denoted by $\frak P$. For each
$\pi\in\frak P$, choose an element $g^\pi=(g^\pi_\bi)_\bi\in G$ such that its
conjugacy class corresponds to $\pi$. Thus, we have
$$|G_{g^\pi}|=\prod_{\bi\in I} |GL_{a_\bi}(\fqi)_{g^\pi_\bi}|
=\prod_{\bi\in I}\prod_{\vphi\in\Phi(q^{\ez_\bi})}
q^{\ez_\bi d(\vphi)\lr{\pi^\bi(\vphi),\,\pi^\bi(\vphi)}}
b_{\pi^\bi(\vphi)}(q^{-\ez_\bi d(\vphi)}).$$
Clearly, for each $\bi\in I$, $\pi^\bi$ corresponds to the conjugacy class of
$g_\bi^\pi$ in $GL_{a_\bi}(\fqi)$. Let $0\leq s\leq\ez_\bi-1$
and denote by $\pi^\bi[s]$ the partition function $\Phi(q^{\ez_\bi})\ra\cp$
corresponding to the conjugacy class of $(g_\bi^\pi)^{[s]}$.
Then we get
$$|X^{g^\pi}|=\prod_{\brho:\bi\ra\bj}|X^{g^\pi}_\brho|
=\prod_{\brho:\bi\ra\bj}q^{\ez_\brho\sum_{\vphi\in\Phi(q^{\ez_\bi}),\,
\psi\in\Phi(q^{\ez_\bj})}d(\vphi,\psi)\lr{\pi^\bi[r_\brho](\vphi),\,
\pi^\bj[s_\brho](\psi)}}.$$
Finally, we deduce
\begin{equation}\label{number}
M_{Q,\sz}(\az,q)=\sum_{\pi\in\frak P}\frac{|X^{g^\pi}|}{|G_{g^\pi}|}
=\sum_{\pi\in\frak P}\frac
{\prod_{\brho:\bi\ra\bj}q^{\ez_\brho\sum_{\vphi\in\Phi(q^{\ez_\bi}),\,
\psi\in\Phi(q^{\ez_\bj})}d(\vphi,\psi)
\lr{\pi^\bi[r_\brho](\vphi),\,\pi^\bj[s_\brho](\psi)}}}
{\prod_{\bi\in I}\prod_{\vphi\in\Phi(q^{\ez_\bi})}
q^{\ez_\bi d(\vphi)\lr{\pi^\bi(\vphi),\,\pi^\bi(\vphi)}}
b_{\pi^\bi(\vphi)}(q^{-\ez_\bi d(\vphi)})}.
\end{equation}

An orientation of the underlying graph $\ol Q$ of $Q$ is called
$\sz$-admissible\footnote{The $\sz$-admissible orientations of the
underlying graph of $Q$ correspond to the orientations of the
underlying graph of $\ggz$.} if it is compatible with the graph
automorphism of $\ol Q$ induced by $\sz$. Obviously, the
orientation of $Q$ itself is $\sz$-admissible.

\begin{thm} The number $M_{Q,\sz}(\az,q)$ is a polynomial in $q$ with rational
coefficients and independent of the $\sz$-admissible orientation of $Q$.
\end{thm}

\begin{pf} By (\ref{number}), $M_{Q,\sz}(\az,q)$ is a rational function in $q$
over $\bbq$. Since $M_{Q,\sz}(\az,q)\in\bbz$ for all power $q=p^s$ of a prime $p$,
it follows that $M_{Q,\sz}(\az,q)$ is a polynomial. Further, whenever $i_0\in\bi$ and
$\rho_0\in \brho$ for $\bi\in I$ and $\brho\in\ggz_1$ are fixed,
the numbers $r_\brho$ and $s_\brho$ for each $\brho\in\ggz_1$ are clearly independent
of the orientation of $\brho$, we have again by (\ref{number}) that
$M_{Q,\sz}(\az,q)$ is independent of the orientation of $\ggz$, i.e., the
$\sz$-admissible orientation of $Q$.
\end{pf}

By induction on the height $\text{ht}\az:=\sum_{\bi\in I}a_\bi$ of
$\az=\sum_{\bi\in I}a_\bi\bi\in\dt(\ggz)^+$, we deduce the
following (see \cite{Hua2} for the case of natural
$\bbf_q$-modulated quivers).

\begin{cor} The number $I_{Q,\sz}(\az,q)$ of isoclasses of indecomposable
$A^F$-modules of dimension vector $\az$ is a polynomial in $q$ with
rational coefficients and independent of the $\sz$-admissible orientation
of $Q$.
\end{cor}

\section{Roots and indecomposable $F$-stable representations}

In this last section, we present an application to Lie theory. We
keep the notation introduced in \S9. Thus, $(Q,\sz)$ denotes an
ad-quiver and $\fM_{Q,\si}$ is the associated modulated quiver
with the underlying valued quiver $\ggz=\ggz(Q,\sz)=(I,\ggz_1)$.

The quiver $Q$ defines a symmetric generalized Cartan matrix
$C_Q=(a_{ij})_{i,j\in Q_0}$ given by
$$a_{ij}=\left\{ \begin{array}{ll}    2 \;\; &\mbox{if}\;\; i=j \\
        -|\{\mbox{arrows between $i$ and $j$}\}| \;\; &\mbox{if}\;\; i\ne j
        \end{array}\right.$$
while the valued quiver $\ggz$ defines a symmetrizable generalized
Cartan matrix $C_\Ga=(b_{\bi\bj})_{\bi,\bj\in I}$ given by
$$b_{\bi\bj}=\left\{ \begin{array}{ll}
2 \;\; &\mbox{if}\;\; \bi=\bj \\
-\sum_{\brho}\ez_\brho/\ez_\bi \;\; &\mbox{if}\;\; \bi\ne \bj
                 \end{array}\right.$$
where the sum is taken over all arrows $\brho$ between $\bi$ and
$\bj$ (see (\ref{GCM})). In fact, all symmetrizable generalized
Cartan matrix can be obtained in this way.

Let  $\dt(\ggz)\subset \bbz I$ be the root system associated with
the valued quiver $\ggz$, or equivalently, the root system of
Kac-Moody algebra associated with the Cartan matrix $C_\Gamma$
(see \cite{K1} and \cite{K3} for definition). We shall write
$\dt(Q)$ for $\dt(\Gamma)$ if $\si=1$ and $\dt(\ggz)^+$ for the
positive subsystem. For $\bz\in\dt(Q)^+$, let $t\geq 1$ be the
minimal number satisfying $\sz^t(\bz)=\bz$. We call $t$ the
$\si$-{\it period} of $\bz$, denoted $p(\bz)=p_\si(\bz)$.

Recall from (\ref{hatsi}) the group isomorphism $\hsz:(\bbz
Q_0)^\sz\ra \bbz I$ which induces a map from $\dt(Q)$ to
$\dt(\Ga)$ in the following.


\begin{lem} \label{roots} (\cite{T},\cite[Prop.4]{Hub1}) Let $\bz\in\dt(Q)$
and let
$${\tilde \bz}:=\bz+\sz(\bz)+\cdots+\sz^{t-1}(\bz)\in(\bbz Q_0)^\sz,$$
where $t=p(\bz)$. Then $\bz\lmto\hsz(\tilde \bz)$ induces a
surjective map $\dt(Q)\ra\dt(\ggz)$. Moreover, if $\hsz(\tilde
\bz)$ is real, then $\bz$ is real and is unique up to $\sz$-orbit.
\end{lem}


Let $A$ be the path algebra of $Q$ and $F$ the Frobenius morphism
induced from $\si$.

\begin{prop} \label{Kac} The dimension vector of each indecomposable
$A^F$-module lies in $\dt(\ggz)^+$. Moreover, if $\az\in\dt(\ggz)^+$
is real, there is a unique, up to isomorphism, indecomposable $A^F$-module
with dimension vector $\az$. In other words, we have
\begin{itemize}
\item[(a)] If $I_{Q,\si}(q,\az)\neq0$, then $\az\in\dt(\ggz)^+$.
\item[(b)] If $\az\in\dt(\ggz)^+$ is real, then $I_{Q,\si}(q,\az)=1$.
\end{itemize}
\end{prop}

\begin{pf} Let $X$ be an indecomposable $A^F$-module
with $\End_{A^F}(X)/\rad(\End_{A^F}(X))\cong \bbf_{q^r}$. By
Theorem \ref{INDF}, there is an indecomposable $A$-module $M$ such
that $M=M^{[r]}$ and
$$\bX:=X\otimes_{\bbf_q}k\cong M\oplus M^{[1]}\oplus\cdots M^{[r-1]}.$$
Moreover, $M, M^{[1]}, \cdots$, and $M^{[r-1]}$ are pairwise
non-isomorphic. By Kac's theorem, the dimension vector $\bz:=\udim
M$ lies in $\dt(Q)^+$. Let $t$ be minimal such that
$\sz^t(\bz)=\bz$. Since $\udim M^{[1]}=\sz(\udim M)$, we have
$t|r$ and
$$\udim \bX=\bz+\sz(\bz)\cdots+\sz^{r-1}(\bz)=\frac{r}{t}{\tilde \bz}.$$
If $\bz$ is an imaginary root, then $\hsz(\tilde \bz)$ is an
imaginary root in $\dt(\ggz)^+$. We conclude that
$$\udim X=\hsz(\udim \bX)=\frac{r}{t}\hsz(\tilde \bz)$$
is an imaginary root in $\dt(\ggz)^+$. If $\bz$ is real, then
$r=t$. This implies that $\udim X=\hsz(\tilde \bz)$ is a real root
in $\dt(\ggz)^+$.

Now let $\az=\hsz(\tilde \bz)$ be a real root. Then, by Lemma
\ref{roots}, $\bz$ is real. Again by Kac's theorem, there is a
unique indecomposable $A$-module $M$ with dimension vector $\bz$.
Since $\udim M^{[t]}=\sz^t(\udim M)=\udim M$, where $t$ is minimal
such that $\sz^t(\bz)=\bz$, we have $M^{[t]}\cong M$. Hence, by
Theorem \ref{INDF}, we obtain an $F$-stable module
$${\tilde M}\cong M\oplus M^{[1]}\oplus\cdots\oplus M^{[t-1]}$$
such that ${\tilde M}^F$ is indecomposable $A^F$-module whose
dimension vector is $\hsz(\tilde \bz)=\az$. The uniqueness of such
an indecomposable module follows from the fact that $\bz$ is
unique up to $\sz$-orbit.
\end{pf}

\begin{rems}\label{Kconj} (1) The statements (a) and (b) in \ref{Kac} together with

\vspace{.3cm}
(c) If $\az\in\dt(\ggz)^+$ is imaginary, then the {\it polynomial
} $I_{Q,\si}(q,\az)\neq0$. \vspace{.3cm}

\noindent constituent the so-called Kac's Theorem (see \cite{K1})
when $\si=1$.  The proposition partially generalizes Kac's theorem
for quivers (over a finite field) to a result for modulated
quivers defined by $(Q,\si)$ with $\si\neq1$. The generalization
of Kac's Theorem was first formulated by Hua in \cite[Thm
4.1]{Hua2}. Hubery in \cite{Hub1} provides a proof by using a
classification of the so-called ii-indecomposable representations
of an ad-quiver. However, both \cite{Hua2} and \cite{Hub1} treat
{\it only} natural $\bbf_q$-modulated quivers (see Remark
\ref{Hua-Hub}).

(2) By using Ringel-Hall algebra approach, it has been proved in
\cite{DX} that, for any prime power $q$ and dimension vector
$\az$, the {\it number} $I_{Q,\sz}(\az,q)\not=0$ if and only if
$\az\in\dt(\ggz)^+$. This is stronger than the statement (c) for
imaginary roots. It should be still interesting to find a direct
proof for (c) in the modulated quiver case.
\end{rems}

We end the paper with the following conjecture which provides a
direct link of Kac's theorem with its generalization.

\begin{con} Let $A$ be the path algebra of an
ad-quiver $(Q,\si)$ without oriented cycles. Let $\bz\in\De(Q)^+$.
Then there exists an indecomposable $A$-module $M$ with dimension
vector $\bz$ such that $p_F(M)=p_\si(\bz)$.
\end{con}

We claim that this conjecture implies \ref{Kconj}(c) immediately.
Indeed, for any $\al\in\De(\Ga)^+$, there exists a
$\bz\in\De(Q)^+$ such that $\al=\hat\si(\tilde\bz)$. By Kac's
theorem, we have $I_{Q,1}(q,\bz)\neq0$ for some $q$. Now the
existence of such an $M$ guarantees that the associated $F$-stbale
module $\tilde M$ defined in (\ref{Mtil}) has dimension vector
$\tilde \bz$. Thus, the indecomposable $\tilde M^F$ has dimension
vector $\al$. Therefore, $I_{Q,\si}(q,\al)\neq0$.


\begin{thebibliography}{9}


\bibitem{ARS} M. Auslander, I. Reiten, and S.O. Smal$\phi$,
{\em Representation Theory of Artin Algebras}, Cambridge Studies in Advanced
Mathematics: {\bf 36}. Cambridge University Press, Cambridge, 1995.

\bibitem{Be} D. Benson, {\em Representations and Cohomology, Vol
I}, Cambridge Studies in Advanced Mathematics: {\bf 30}. Cambridge
University Press, Cambridge, 1995.


\bibitem{CR} C.W. Curtis and I. Reiner, {\em Methods of representation theory.
With applications to finite groups and orders}, Vol. I.  A Wiley-Interscience
Publication, New York, 1981.

\bibitem{DX} B. Deng and J. Xiao,
{\em A new approach to Kac's theorem on representations of valued
quivers}, to appear in Math. Z.

\bibitem{DD1} B. Deng and J. Du, {\em Monomial basse for quantum
affine ${\frak {sl}}_n$}, to appear.

\bibitem{DD2} B. Deng and J. Du, {\em On bases of quantized enveloping algebras}, to appear.

\bibitem{DM} F. Digne and J. Michel,
{\em Representations of finite groups of Lie type}, London Math. Soc. Student Texts,
21. Cambridge University Press, Cambridge, 1991. .

\bibitem{DR1} V. Dlab and C.M. Ringel, {\em On
algebras of finite representation type}, J. Algebra {\bf 33},
(1975), 306-394.

\bibitem{DR2} V. Dlab and C.M. Ringel, {\em Indecomposable
representations of graphs and algebras}, Memoirs Amer. Math. Soc. {\bf 6}
no. 173, 1976.

\bibitem{DF} P.W. Donovan and M.R. Freislich,
{\em The representation theory of finite graphs and associated algebras},
Carleton Math. Lecture Notes {\bf 5}, 1973.

\bibitem{Ga} P. Gabriel,
{\em Unzerlegbare Darstellungen I}, Manuscripta Math. {\bf 6}
(1972), 71-103.

\bibitem{Ga2} P. Gabriel,
{\em Auslander-Reiten sequences and representation-finite algebras},
Lecture Notes in Math., {\bf 831}, Springer-Verlag, Berlin-Heidelberg-New
York (1980), 1-71.

\bibitem{GR} P. Gabriel and A. V.Roiter, {\em Representations of finite dimensional algebras},
Spriner-Verlag, Berlin Heidelberg, 1997.

\bibitem{Hua1} J. Hua, {\em Representations of quivers over finite
fields}, Ph.D. thesis, University of New South Wales, 1998.

\bibitem{Hua2} J. Hua, {\em Numbers of representations of valued quivers over
finite fields}, preprint, Universit\"at Bielefeld, 2000
(www.mathematik.uni-bielefeld.de/$^\sim$sfb11/vquiver.ps).

\bibitem{HL} J. Hua and Z. Lin,
{\em Generalized Weyl denominator formula}, In: Representations
and quantizations, Proceedings of the International Conf. on
Representation Theory (Sahnghai, 1998), 247-261, China High. Educ.
Press, Beijing, 2000.

\bibitem{Hub1} A. Hubery, {\em Quiver representations respecting a quiver
automorphism: a generalisation of a theorem of Kac}, preprint, 2002
(math.RT/0203195).

\bibitem{Hub2} A. Hubery, {\em Representations of quivers respecting
a quiver automorphism and a theorem of Kac}, Ph.D. thesis, University
of Leeds, August 2002.

\bibitem{K1} V. Kac, {\em Infinite root systems, representations of graphs
and invariant theory}, Invent. Math. {\bf 56}(1980), 57-92.

\bibitem{K2} V. Kac, {\em Root systems, representations of quivers and
invariant theory}, Lecture Notes in Mathematics {\bf 996},
Springer-Verlag, 1982, 74-108.

\bibitem{K3} V. Kac, {\em Infinite dimensional Lie algebras}, Third edition,
Cambridge University Press, 1990.

\bibitem{M} I.G. Macdonald,
{\em Symmetric functions and Hall polynomials}, 2nd edition, Clarendon Press,
Oxford, 1995.

\bibitem{L93} G. Lusztig, {\em Introduction to quantum groups}, Progress in
Math. {\bf 110}, Birkh\"auser, 1993.

\bibitem{L98} G. Lusztig, {\em Canonical bases and Hall algebras,}  Representation
theories and algebraic geometry,  365-399, Kluwer Acad. Publ., Dordrecht, 1998.

\bibitem{N} L.A. Nazarova,
{\em Representations of quivers of infinite type},
Math. USSR Izvestija Ser. Mat. {\bf 7} (1973), 752-791.

\bibitem{Re} M. Reineke,
{\em The quantic monoids and degenerate quantized enveloping
algebras}, preprint, 2002.

\bibitem{T} T. Tanisaki,
{\em Foldings of root systems and Gabriel's theorem},  Tsukuba J. Math.
{\bf 4}(1980), 89-97.


\end{thebibliography}
\end{document}